\documentclass[11pt,reqno,twoside]{article}
\usepackage[british]{babel}

\usepackage{a4wide,amscd,amsmath,amsthm,amssymb}
\usepackage{mathrsfs}
\usepackage[latin1]{inputenc}
\usepackage{euscript}

\DeclareFontFamily{OT1}{rsfs}{}
\DeclareFontShape{OT1}{rsfs}{n}{it}{<-> rsfs10}{}
\DeclareMathAlphabet{\curly}{OT1}{rsfs}{n}{it}

\DeclareMathOperator{\Ker}{Ker}
\DeclareMathOperator{\Spin}{Spin}
\DeclareMathOperator{\hcf}{gcd}
\DeclareMathOperator{\pr}{pr}
\renewcommand{\Re}{\operatorname{Re}}
\renewcommand{\Im}{\operatorname{Im}}

\newcommand{\RE}{\mathbb{R}}
\newcommand{\CX}{\mathbb{C}}
\newcommand{\ZE}{\mathbb{Z}}
\newcommand{\QU}{\mathbb{Q}}
\newcommand{\CP}{\mathbb{C}P}
\newcommand{\CPA}{\mathbb{C}P^n_a}
\newcommand{\CPW}{\mathbb{C}P^5_{1,1,1,1,4,4}}
\newcommand{\OH}{\curly{O}}
\newcommand{\Ga}{\Gamma}
\newcommand{\tgamma}{\tilde{\gamma}}
\newcommand{\eps}{\varepsilon}
\newcommand{\p}{\partial}
\newcommand{\we}{\wedge}
\newcommand{\AM}{\EuScript{A}}
\newcommand{\J}{\mathcal{J}}
\newcommand{\HH}{\mathcal{H}}
\newcommand{\minus}{\smallsetminus}
\newcommand{\cupprod}{\mathbin{\raisebox{-2pt}{$\smallsmile$}}}

\newcommand{\eus}{\EuScript}

\renewcommand{\phi}{\varphi}
\let\emptyset\varnothing

\newcommand{\tPhi}{\tilde\Phi}
\newcommand{\te}{\tilde\eta_3}
\newcommand{\trho}{\tilde{\rho}}
\newcommand{\tV}{\tilde{V}}
\newcommand{\tD}{\tilde{D}}
\newcommand{\oM}{\overline{M}}
\newcommand{\oW}{\overline{W}}

\newtheorem{thrm}{Theorem}[section]
\newtheorem{prop}[thrm]{Proposition}

\newtheorem{cor}[thrm]{Corollary}
\theoremstyle{definition}
\newtheorem{defi}[thrm]{Definition}
\theoremstyle{remark}
\newtheorem{rmk}[thrm]{Remark}

\begin{document}

\title{Asymptotically cylindrical manifolds with holonomy $\Spin(7)$. I}
\author{Alexei Kovalev
\\[5pt]
DPMMS, University of Cambridge,\\
   Centre for \mbox{Mathematical} Sciences,\\
   \mbox{Wilberforce} Road, Cambridge CB3 0WB, UK
\\[5pt]
\hbox{\tt a.g.kovalev@dpmms.cam.ac.uk}
}
\date{\null}
\maketitle

\begin{abstract}
We construct examples of asymptotically cylindrical Riemannian 8-manifolds
with holonomy group $\Spin(7)$. To our knowledge, these are the first such
examples. The construction uses an extension to asymptotically cylindrical
setting of Joyce's existence result for torsion-free $\Spin(7)$-structures.
One source of examples arises from `Fano-type' K\"ahler 4-orbifolds with
smooth anticanonical Calabi--Yau 3-fold divisors and with compatible
antiholomorphic involution. We give examples using weighted projective
spaces and calculate basic topological invariants of the resulting
$\Spin(7)$-manifolds.
\end{abstract}

\markboth{\hfil\mbox{\sc alexei kovalev}\hfil\hfil}
{\hfil\hfil\mbox{\sc asymptotically cylindrical spin(7)-manifolds}\hfil}

The holonomy group $\Spin(7)$ occurs in 8 dimensions as one of the two
special cases in Berger's classification of the holonomy of the
Levi--Civita connections; the other being the holonomy group $G_2$ in
7~dimensions. Both special holonomy groups $\Spin(7)$ and $G_2$ are
naturally related to the algebra of octonions, or Cayley numbers; the group
$\Spin(7)$ arises as the stabilizer of the triple cross
product on $\RE^8$ \cite{HL,salamon-walpuski}. One the other hand, the
holonomy representation of $\Spin(7)$ is isomorphic to the spin
representation.

The metrics with holonomy $\Spin(7)$ and $G_2$ are Ricci-flat and are the
only holonomy groups in Berger's list corresponding to Ricci-flat metrics
that are {\em not} K\"ahler. More explicitly, the other Ricci-flat
holonomy groups are $SU(n)$ and $Sp(n)\subset SU(2n)$; these induce a
parallel, integrable complex structure on a base manifold. The proper
inclusion of Riemannian holonomy groups $SU(4)\subset\Spin(7)$ accounting
for holonomy reduction with no compatible complex structure is unique to
real dimension~8.

The first examples of Riemannian manifolds with holonomy $\Spin(7)$ were
obtained by Bryant~\cite{bryant}. Complete examples were subsequently
constructed by Bryant and Salamon~\cite{bryant-salamon} and by Gibbons, Page
and Pope~\cite{GPP}. The first compact manifolds with holonomy $\Spin(7)$
are due to Joyce \cite{joyce-inv,joyce-jdg}.

In this paper, we construct examples of complete asymptotically
cylindrical 8-manifolds with holonomy equal to~$\Spin(7)$. To the author's
knowledge, these are the first such examples. The previously known complete
non-compact holonomy $\Spin(7)$-metrics are of asymptotically conical type
with a different (non-linear polynomial) volume growth.

Examples of asymptotically cylindrical Calabi--Yau manifolds
with holo\-nomy $SU(n)$ were constructed in \cite{g2paper}, \cite{KL}
and~\cite{CHNP}. A construction of asymptotically cylindrical manifolds
with holonomy $G_2$ was given in~\cite{KN}.

The metrics with holonomy $\Spin(7)$ can be determined via a differential
self-dual 4-form $\Phi$ point-wise modelled on the `Cayley form' on $\RE^8$.
Any such 4-form $\Phi$ induces a $\Spin(7)$-structure on an 8-manifold,
hence a metric as $\Spin(7)$ is a subgroup of $SO(8)$. The holonomy of this
metric is contained in $\Spin(7)$ precisely when the $\Spin(7)$-structure
is torsion-free, a condition equivalent to $\Phi$ being a closed 4-form.
In practical examples, the task of constructing a suitable 4-form leads
to a non-linear PDE. A criterion for an asymptotically cylindrical
torsion-free $\Spin(7)$-structure to have the full holonomy $\Spin(7)$ was
proved by Nordstr\"om \cite{jn-thesis}; for simply-connected examples
constructed here it amounts to the vanishing of the first Betti number of
the cross-section of cylindrical end.

The asymptotically cylindrical $\Spin(7)$-manifolds in this paper are
obtained by modifying Joyce's construction of compact $\Spin(7)$-manifolds
in~\cite{joyce-jdg}. Joyce uses in {\it op.cit.} quotients of
Calabi--Yau complex 4-dimensional orbifolds with suitable isolated
singularities by anti-holomorphic involution and resolves the singularities
to obtain smooth 8-manifolds with $\Spin(7)$-structures having small
torsion.
If, instead, a K\"ahler 4-orbifold has a smooth Calabi--Yau 3-fold divisor in
the anticanonical class, then the complement of this divisor, under suitable
conditions, admits an asymptotically cylindrical metric with holonomy $SU(4)$
cf.~\cite{g2paper,HHN}. It is also possible to choose a `compatible'
anti-holomorphic involution to obtain asymptotically cylindrical 
$\Spin(7)$-orbifolds. Furthermore, a large part of the resolution of
singularities and the existence argument for torsion-free
$\Spin(7)$-structures in \cite[Chap. 13 and 15]{joyce} is carried out in way
which can be extended to asymptotically cylindrical setting with little
further work. The situation here has some analogy with \cite[\S 3]{KN},
where a generalization similar type is provided to obtain an existence
result for asymptotically cylindrical $G_2$-structures. In particular, the
main additional argument required in both the $\Spin(7)$ and $G_2$ cases is
to show an exponential rate of convergence to cylindrical asymptotic
structure. The task of applying the existence result in examples then
amounts to finding suitable orbifolds (see Definition~\ref{config}) and
this can be treated as a question in complex algebraic geometry.

The paper is organized as follows. We begin with introducing the necessary
background results on torsion-free $\Spin(7)$-structures and holonomy in
\S\ref{synopsis} and on asymptotically cylindrical manifolds including
aspects of Hodge theory in \S\ref{as.cyl}. Then, in \S\ref{proper.spin7},
we explain that the holonomy group of a simply-connected asymptotically
cylindrical manifold is topologically determined. This includes a criterion
for full holonomy $\Spin(7)$ and a description of the moduli space of
asymptotically cylindrical $\Spin(7)$-manifolds, the results due to
Nordstr\"om~\cite{jn-thesis}. We also give an interpretation of some
results of~\cite{HHN} in the context of $\Spin(7)$-structures.
An asymptotically cylindrical extension of Joyce's existence theorem for
$\Spin(7)$ structures is proved in~\S\ref{exist}.
The class of K\"ahler 4-orbifolds required for making asymptotically
cylindrical $\Spin(7)$ manifolds is introduced in \S\ref{orbifolds} 
which also contains justification of the construction.
Finally, in \S\ref{examples}, we show how to implement the construction. We
restrict attention to simple examples and calculate their basic topological
invariants and the dimension of respective moduli space.

The asymptotically cylindrical $\Spin(7)$-manifolds in this paper have
Betti numbers $b^1=b^2=b^3=0$ and with their cross-sections having
$b^1=b^2=0$. This is in fact the maximal set of vanishing Betti numbers
possible for asymptotically cylindrical $\Spin(7)$-manifolds, whether or
not the holonomy group is all of $\Spin(7)$, as the cohomology class of
$\Spin(7)$-form is never zero and induces non-trivial harmonic 3-form
defining a $G_2$-structure on the (compact) cross-section `at infinity'.
On the other hand, the necessary condition for $\Spin(7)$ holonomy for
8-manifolds of this type is the vanishing of $b^1$ for both the 8-manifold
and its cross-section. Further, the holonomy of the cross-section induced
by the $G_2$-structure at infinity may be exactly $G_2$ or a proper
subgroup, depending on whether or not the fundamental group is finite.
It is infinite in the present examples and the holonomy the $G_2$-manifold
at infinity is $\ZE_2\ltimes SU(3)$.

The author has a method which produces examples of asymptotically cylindrical
holonomy $\Spin(7)$ manifolds where the $G_2$-manifold at infinity is
irreducible with full holonomy~$G_2$. These developments and some
applications are treated in the companion paper~\cite{II}.

\section{$\Spin(7)$-structures on 8-manifolds}
\label{synopsis}

We give a short summary of some background results on the Riemannian geometry
in dimension 8 associated with the structure group $\Spin(7)$.
For more details, see \cite{HL,joyce,salamon-walpuski,salamon}.

The key role in defining the $\Spin(7)$ holonomy is played by a particular
4-form, sometimes called the Cayley 4-form, on the Euclidean $\RE^8$
\begin{multline}\label{std.spin7}
\Phi_0=dx_{1234}+dx_{1256}+dx_{1278}+dx_{1357}-dx_{1368}
-dx_{1458}-dx_{1467}\\
-dx_{2358}-dx_{2367}-dx_{2457}+dx_{2468}+dx_{3456}+dx_{3478}+dx_{5678}.
\end{multline}
The form $\Phi_0$ arises by considering $\RE^8$ as the algebra of
octonions, or Cayley numbers, and setting
$\Phi_0(x,y,z,w)=\frac12\langle x(\bar{y}z)-z(\bar{y}x),w\rangle$.
The stabilizer of $\Phi_0$ in the $GL(8,\RE)$ action
is a subgroup of $SO(8)$ isomorphic to $\Spin(7)$
\cite[Prop. IV.1.36]{HL}. The Cayley form is
self-dual $*\Phi_0=\Phi_0$ with respect to the Euclidean metric.

For an oriented 8-manifold~$M$, define a subbundle of {\em admissible} 4-forms
$\AM M\subset \Lambda^4 T^* M$ 
with the fibre $\AM_p M$ at each $p\in M$ consisting of 4-forms that can be
identified with $\Phi_0$ via an orientation-preserving isomorphism
$T_pM\to\RE^8$. The fibres of $\AM M$ are diffeomorphic to the orbit
$GL_+(8,\RE)/\Spin(7)$ of $\Phi_0$, a 43-dimensional submanifold of
the 70-dimensional vector space $\Lambda^4 (\RE^8)^*$.

Every admissible 4-form $\Phi\in\Gamma(\AM M)$ defines a $\Spin(7)$-structure
on $M$, hence a metric $g=g(\Phi)$ and orientation and a Hodge star $*_\Phi$,
so that ${*_\Phi}\Phi=\Phi$. We shall sometimes slightly inaccurately
say that that $\Phi$ {\em is} a $\Spin(7)$-structure.

Here are two foundational results on the $\Spin(7)$-structures, the first is
taken from \cite[Lemma~12.4]{salamon} and the second from~\cite{bonan}.

\begin{thrm}
Let $M$ be an 8-manifold and $\Phi\in\Ga(\AM M)$ a $\Spin(7)$-structure on~$M$.
Then

(i)
the holonomy of $g(\Phi)$ is contained in $\Spin(7)$ if and only if $d\Phi=0$,
and

(ii)
if the holonomy of a metric $g$ on~$M$ is contained in $\Spin(7)$, then
$g$ is Ricci-flat.
\end{thrm}

The $\Spin(7)$-structure induced by a closed admissible 4-form $\Phi$ is
said to be {\em torsion-free}; the condition $d\Phi=0$ in this case is
equivalent to vanishing of the intrinsic torsion. We shall call an
8-manifold $(M,\Phi)$ endowed with a torsion-free $\Spin(7)$-structure
a {\em $\Spin(7)$-manifold}.

Each $\Spin(7)$-structure on an 8-manifold $M$ defines point-wise orthogonal
decompositions of the bundles of differential $r$-forms corresponding to
irreducible components of the induced representation of $\Spin(7)$ on
$\Lambda^r (\RE^8)^*$ (\cite{bryant} or \cite{fernandez}).
In particular,
\begin{subequations}\label{rep4}
\begin{align}
\Lambda^2 T^* M &= \Lambda^2_7 T^* M \oplus \Lambda^2_{21} T^* M,
\label{2sub}\\
\Lambda^3 T^* M &= \Lambda^3_8 T^* M \oplus \Lambda^3_{48} T^* M,
\\
\Lambda^4 T^* M &= \Lambda^4_- T^* M \oplus \Lambda^4_+ T^* M,
\qquad
\Lambda^4_- T^* M = \Lambda^4_{35} T^* M,
\\
\Lambda^4_+ T^* M &=
\Lambda^4_1 T^* M \oplus \Lambda^4_7 T^* M \oplus \Lambda^4_{27} T^* M
\end{align}
\end{subequations}
where the lower indices indicate ranks of subbundles and the fibres of
$\Lambda^4_\pm$ are the $\pm 1$-eigenspaces of $*_\Phi$.
We may sometimes abbreviate $\Lambda^r_{j} T^* M$ as $\Lambda^r_{j}$ if no
confusion is likely.

The subbundles of 2-forms in~\eqref{2sub} may be determined by
\begin{equation}\label{sub2forms}
\Lambda^2_7 = \{\alpha : *(\Phi\we\alpha)=3\alpha\},\qquad
\Lambda^2_{21} = \{\alpha : *(\Phi\we\alpha)=-\alpha\}.
\end{equation}
The fibres of $\Lambda^4_1$ are spanned by $\Phi$ and $\Lambda^3_8$ is the
bundle of 3-forms $v\lrcorner\Phi$ for $v\in TM$. The inner product
$g(\Phi)$ is determined by
$$
\langle u,v\rangle\, \Phi\we\Phi =
2\, (u\lrcorner\Phi)\we *_\Phi(v\lrcorner\Phi).
$$

For each $\Phi_p\in \AM_p M$, $p\in M$, the tangent space to the orbit
$(\AM M)_p\subset \Lambda^4T^*_p M$ at $\Phi_p$ is 
$(\Lambda^4_1 \oplus \Lambda^4_7 \oplus \Lambda^4_{35})_p$
and the normal space is the orthogonal complement $(\Lambda^4_{27})_p$, where
the subspaces $(\Lambda^4_k)_p$ are determined by the $\Spin(7)$
structure~$\Phi_p$. Using this, one can construct a tubular open
$\eps$-neighbourhood $\eus{T}M$ of $\AM M$ in $\Lambda^4 T^* M$ consisting
of 4-forms $\chi_p=\Phi_p+\psi_p$ such that $\Phi_p\in\AM_p M$,
$\psi_p\in(\Lambda^4_{27})_p$ and $|\psi_p|<\eps$ (with
$(\Lambda^4_{27})_p$ and the point-wise norm induced again by $\Phi_p$).  
Furthermore, for each $\chi_p\in\eus{T}M$, $\Phi_p$ and $\psi_p$ are
uniquely determined by the latter conditions, if $\eps>0$ was chosen
sufficiently small. The $\eus{T}M$ is a fibre bundle of 4-forms (it is not
a vector bundle) and there is a well-defined smooth projection
$$
\Theta:\chi_p=\Phi_p+\psi_p\in\eus{T}M\to \Phi_p\in\AM M
$$
The definition of $\Theta$ does not use any `background'
$\Spin(7)$-structure on~$M$.

We note for later use that the map $\Theta$ has a local expansion at each
$\Spin(7)$-structure $\Phi$,
$$
\Theta(\Phi+\psi)=
\Phi+\pi_1(\psi)+\pi_7(\psi)+\pi_{35}(\psi)-F(\psi),
\qquad
\psi\in\Omega^4(M),\quad \|\psi\|_{C^0}<\eps,
$$
where $\pi_i$ denotes the $\Lambda^4_i$ component of a 4-form in the
decomposition~\eqref{rep4} determined by~$\Phi$ and the norm is computed
using the metric $g(\Phi)$.
The remainder $F$ satisfies quadratic estimates
\begin{align}
|F(\psi')-F(\psi'')| < &\; C_1|\psi'-\psi''|(|\psi'|+|\psi''|),
\label{quadratic}\\
\begin{split}
|\nabla(F(\psi')-F(\psi''))| < 
&\; C_2\bigl(|\psi'-\psi''|(|\psi'|+|\psi''|)|d\Phi|
 +|\nabla(\psi'-\psi'')|(|\psi'|+|\psi''|)
\\
&\; +|\psi'-\psi''|(|\nabla\psi'|+|\nabla\psi''|)\bigr),
\end{split}\label{quadratic2}
\end{align}
with $C_1>0$, $C_2>0$ independent of $\psi',\psi''$ \cite[Prop.~10.5.9]{joyce}.

The holonomy of a $\Spin(7)$-manifold $M$ may be contained in a subgroup
$SU(4)\subset\Spin(7)$. In the notation of~\eqref{std.spin7} the holonomy
group $SU(4)$ is the subgroup of $GL(8,\RE)$ preserving the standard
symplectic form $\omega_0=dx_1\we dx_2+\ldots+dx_7\we dx_8$ and the real
part of the 4-form $\theta_0=(dx_1+idx_2)\we\ldots\we(dx_7+idx_8)$.
An 8-manifold $M$ with holonomy in $SU(4)$ carries a parallel integrable
complex structure with trivial canonical bundle and a Ricci-flat K\"ahler
metric, i.e.\ $M$ is a Calabi--Yau complex 4-fold. The parallel forms
corresponding to $\omega_0$ and $\theta_0$ are the K\"ahler form and a
non-vanishing holomorphic $(4,0)$-form, sometimes called a holomorphic
volume form, normalized so that $3\theta\we\bar\theta=2\omega^4$. These
induce a $\Spin(7)$-structure expressed as
\begin{equation}
\label{su4spin7}
\Phi(\omega,\theta)=\frac12 \omega\we\omega + \Re\theta,
\end{equation}
cf.~\cite[Prop. 13.1.4]{joyce}. A straightforward computation
\begin{equation}\label{kaehler}
*(\omega\we\Phi)=\frac12 *\omega^3= 3\omega,
\end{equation}
shows that the K\"ahler form is a section of the subbundle~$\Lambda^2_7$
determined by $\Phi(\omega,\theta)$.

It will be useful to note one more instance when the holonomy of a
$\Spin(7)$-manifold $M$ reduces to a proper subgroup of~$\Spin(7)$.
The Cayley 4-form on $\RE^8$ may be written as
$\Phi_0=dx_1\we\phi_0 + *_7\phi_0$,
where $\phi_0\in\Lambda^3(\RE^7)^*$ and $*_7$ the Hodge star on the
Euclidean $\RE^7$ with coordinates $x_2,\ldots,x_8$. The stabilizer of
$\phi_0$ in the $GL(7,\RE)$ action is a subgroup of $\Spin(7)$ fixing a 
non-zero vector. This subgroup is isomorphic to the exceptional Lie
group~$G_2$, which also is the stabilizer of the 4-form $*_7\phi_0$.

For a 7-manifold $Y$, a 3-form $\phi\in\Omega^3(Y)$ will be called
{\em stable}, if for each $y\in Y$ there is a linear isomorphism
$j_y:\Lambda^3T^*_y Y\to \RE^7$ smoothly depending on $y$ such that
$j_y^*\phi_0=\phi_y$. Every stable form $\phi$ induces a
$G_2$-structure on~$Y$ and we shall sometimes say that $\phi$ {\em is} a
$G_2$-structure. As $G_2\subset SO(7)$, every $G_2$-structure determines
a metric $g(\phi)$, an orientation and a Hodge star $*_\phi$ on~$Y$.
The holonomy of $g(\phi)$ is contained in $G_2$ if and only if
$$
d\phi=0
\text{ and }
d*_\phi\phi=0.
$$
\cite[Lemma~11.5]{salamon}. In this case, a $G_2$-structure
$\phi$ on $Y$ is said to be torsion-free. A 7-manifold $(Y,\phi)$
endowed with a torsion-free $G_2$-structure will be called a
{\em $G_2$-manifold}.

We shall call a $\Spin(7)$-structure $\Phi_\infty$ on $M=\RE\times Y$
{\em cylindrical} if it induces a product metric $g(\Phi_\infty)=dt^2+g_Y$,
where $t$ is the coordinate on~$\RE$. Equivalently, $\Phi_\infty$ is invariant
under translations in~$t$. Every cylindrical $\Spin(7)$-structure can be
written as
\begin{equation}\label{g2spin7}
\Phi_\infty=dt\we\phi + *_\phi \phi,
\end{equation}
for some $G_2$-structure $\phi$ on $Y$ independent of~$t$,
where $g(\phi)=g_Y$ and $*_\phi$ is the Hodge star of~$g(\phi)$.
A $\Spin(7)$-cylinder $M$ has a parallel vector field
$\frac\p{\p t}$ and the holonomy of $M$ is identified
with the holonomy of the 7-manifold~$Y$. It is not difficult
to check that a cylindrical $\Spin(7)$-structure~\eqref{g2spin7}
on $\RE\times Y$ is torsion-free if and only if the $G_2$-structure $\phi$
on~$Y$ is so. Cf.~\cite[Prop. 13.1.3]{joyce}.

\section{Asymptotically cylindrical metrics}
\label{as.cyl}

Let $M$ be an manifold with cylindrical ends, i.e.\ $M$ is decomposed as
a union of a compact manifold $M_0$ with boundary 
and a half-cylinder $M_\infty=\RE_+\times Y$, with $M_0$ and $M_\infty$
identified along the common boundary $\{0\}\times Y\cong \p M_0$.
The closed manifold $Y$ is called the {\em cross-section} of~$M$.
If $Y$ is connected, we say that $M$ has a single cylindrical end.
We shall use $t$ to denote a smooth function on $M$ which is equal to the
standard $\RE_+$-coordinate on the cylindrical end $M_\infty$ and is
non-positive on~$M_0$.

A $\Spin(7)$-structure $\Phi\in\AM M$ will be called {\em asymptotically
cylindrical with rate~$\lambda$} if there is a cylindrical
$\Spin(7)$-structure $\Phi_\infty$ on $M_\infty$ such that
\begin{equation}\label{acs}
\sup_{\{t\}\times Y}|\nabla_\infty^k (\Phi-\Phi_\infty)|<C_k e^{-\lambda t},
\qquad
t>0,\quad k=0,1,2,\ldots.
\end{equation}
Here the point-wise norms and the covariant derivative are those of
$g_\infty=g(\Phi_\infty)=dt^2+g_Y$. It follows that the metric $g(\Phi)$
approaches a Riemannian product metric $g_\infty$ at an exponential rate,
$$
\sup_{\{t\}\times Y}|\nabla_\infty^k (g-g_\infty)|<C'_k e^{-\lambda t},
\qquad
t>0,\quad k=0,1,2,\ldots ,
$$
so $(M,g(\Phi))$ is an asymptotically cylindrical Riemannian manifold,
with the same rate~$\lambda$. Recall that $\Phi_\infty$ is determined,
via~\eqref{g2spin7}, by a torsion-free $G_2$-structure $\phi$ on~$Y$ and
$g_Y=g(\phi)$.

We shall call an 8-manifold with cylindrical ends endowed with a
torsion-free $\Spin(7)$-structure satisfying~\eqref{acs} an
{\em asymptotically cylindrical $\Spin(7)$-manifold} $(M,\Phi)$ and the
$G_2$-manifold $(Y,\phi)$ the {\em cross-section of $M$ at infinity}.

{\em In this paper, unless stated otherwise, the metric on the
cross-section of any asymptotically cylindrical manifold will
always be understood as the metric at infinity defined above.}

There is a natural class of asymptotically translation-invariant
differential operators defined on manifolds with cylindrical ends,
with coefficients approaching, at an exponential rate, the coefficients of
$t$-independent (translation-invariant) operators on the respective
cylinder. Any differential operator canonically associated to an
asymptotically cylindrical manifold $(M,g)$ will automatically be
asymptotically translation-invariant, in particular the Hodge Laplacian
$\Delta$ has this property.
The respective translation-invariant operator is the Laplacian
$\Delta_\infty$ of the product metric $g_\infty$ on the cylindrical end.
Every differential $r$-form $\alpha$ on $\RE_+\times Y\subset M$ may be
written as a path of forms on~$Y$,
$\alpha=\alpha_r(t) + \alpha_{r-1}(t)\we dt$, and
$$
\Delta_\infty\, \alpha=-\frac{d^2\alpha}{dt^2}
+\Delta_Y\alpha_r+(\Delta_Y\alpha_{r-1})\we dt,
$$
with $\Delta_Y$ the Laplacian on the cross-section.

A general theory of the asymptotically translation-invariant elliptic
operators is developed in \cite{lockhart,LM,MP,melrose}. It includes a
generalization of the Hodge theory isomorphism, between harmonic forms and
the de Rham cohomology classes on closed compact manifolds, to the
asymptotically cylindrical manifolds.

A manifold $M$ with cylindrical ends is non-compact and there are two
standard versions of the de Rham complex, depending on whether one uses
spaces of all differential forms or subspaces of the forms with compact
support. We shall write $H^r(M)$ and $H^r_c(M)$ for the respective
cohomology groups; both are finite-dimensional vector spaces and are
related via the long exact sequence 
\begin{equation}\label{exact}
\ldots\to H^{r-1}(Y)\to H^r_c(M)\to H^r(M)\to H^r(Y)\to\ldots .
\end{equation}
where $Y$, as before, denotes the cross-section of $M$. The Poincar\'e
duality gives isomorphisms $H^r_c(M)\cong H^{\dim M-r}(M)$.

It will be useful to define $H^r_0(M)\subset H^r(M)$ to be the image
of~$H^r_c(M)$ under the `inclusion homomorphism' in~\eqref{exact}, that is,
$H^r_0(M)$ is the subspace of classes represented by compactly supported
closed $r$-forms. We denote $b^r_0(M)=\dim H^r_0(M)$.

Let $\HH^r(M)$ denote the space of all {\em bounded harmonic $r$-forms}
on~$M$.

\begin{thrm}\label{harmonic}
Suppose that $M$ is an asymptotically cylindrical oriented manifold with
rate $\lambda>0$. Let $\eps>0$ be such that $\eps<\lambda$ and $\eps^2$ is
less than any positive eigenvalue of the Hodge Laplacian on differential
forms on~$Y$ of~$M$.
Then

(a)
every $\alpha\in\HH^r(M)$ is smooth and can be asymptotically expressed
on the end of~$M$ as
\begin{equation}\label{asympt}
\alpha|_{M_\infty}=dt\we\alpha_\nu + \alpha_\tau + O(e^{-\eps t}),
\end{equation}
for some harmonic forms $\alpha_\nu\in\HH^{r-1}(Y)$,
$\alpha_\tau\in\HH^r(Y)$ (pulled back via the natural projection) and some
$\eps>0$, with $O(e^{-\eps t})$ understood in the strong sense of
decaying to zero with all derivatives as $t\to\infty$ uniformly on~$Y$.

(b) every $\alpha\in\HH^r(M)$ is closed and co-closed. Define
$$
\HH^r_\tau(M)=\{\alpha\in\HH^r(M):\alpha_\nu=0\}
$$
with $\alpha_\nu=0$ defined by~\eqref{asympt}. Then the map
\begin{equation}\label{hodge}
\alpha\in\HH^r_\tau(M)\to [\alpha]\in H^r(M)
\end{equation}
assigning to a bounded harmonic form its de Rham cohomology class is a
linear isomorphism.

The image under~\eqref{hodge} of the subspace $\HH^r_0(M)\subset\HH^r(M)$
of $L^2$ harmonic forms is $H^r_0(M)$.
\end{thrm}

Part (a) is a direct consequence of the general result in \cite[\S 6]{MP}
on the kernel elements of asymptotically translation-invariant elliptic
operators, cf.\ also \cite[\S 6.4]{melrose}. For (b), see \cite[Prop. 6.14
and 6.18]{melrose}, cf. also \cite[Prop. 4.9]{aps} and
\cite[Theorem 7.9]{lockhart} for the last claim. The results in \cite{melrose}
are given for `exact $b$-metrics' which are equivalent to asymptotically
cylindrical metrics having asymptotic expansions at $t\to\infty$, but the
arguments we require carry over with only cosmetic changes to the present
setting.

\begin{cor}\label{bounded}
Let $M$ be an asymptotically cylindrical oriented manifold with
cross-section~$Y$. Then the dimension of $\HH^r(M)$ is
$b^r(M)+b^r_c(M)-b^r_0(M)$.
\end{cor}
\begin{pf}
Define $\HH^r_\nu(M)=\{\alpha\in\HH^r(M):\alpha_\tau=0\}$
with $\alpha_\tau$ defined in~\eqref{hodge}. We deduce from~\eqref{hodge}
that $*\,\HH^r_\nu(M)=\HH^{\dim M-r}_\tau(M)$, then
$\dim\HH^r_\nu(M)=b^{\dim M - r}(M)=b^r_c(M)$ from Theorem~\ref{harmonic}(b)
and the Poincar\'e duality.

For $\alpha\in\HH^r(M)$, the image of $[\alpha]$ under homomorphism
$\iota^*:H^r(M)\to H^r(Y)$ in~\eqref{exact} is $[\alpha_\tau]$,
from~\eqref{asympt} and the Hodge theory on compact~$Y$. Then from
Theorem~\ref{harmonic}(b) and the exactness of~\eqref{exact} we find that
there exists $\alpha'\in\HH^r_\tau$ with
$\iota^*[\alpha']=[\alpha_\tau]=\iota^*[\alpha]$, so
$\alpha-\alpha'\in\HH^r_\nu(M)$. Thus $\HH^r(M)=\HH^r_\tau(M)+\HH^r_\nu(M)$.
The intersection $\HH^r_\tau(M)\cap\HH^r_\nu(M)=\HH^r_0(M)$ has dimension
$b^r_0(M)$ and the result follows by linear algebra.
\end{pf}

There is a version of the Hodge decomposition theorem for
asymptotically cylindrical manifolds.
\begin{prop}[{cf.~\cite[p.~328]{jn-math.proc}}]\label{decomp}
Let $M$, $\lambda$ and $\eps$ be as in Theorem~\ref{harmonic}. Then for
each $k=0,1,2,\ldots$, every differential form $\xi\in C^k(\Omega^r(M))$
satisfying the exponential decay condition
\begin{equation}\label{k.decay}
\sup_{\{t\}\times Y}|\nabla^j\xi|<C_je^{-\eps t}
\text{ for } t>o,\quad j=0,\ldots,k,
\end{equation}
can be written as
$$
\xi=d\alpha+d^*\beta+\gamma
$$
with $\Delta\gamma=0$, for some unique $d\alpha$, $d^*\beta$, $\gamma$
satisfying~\eqref{k.decay} and pair-wise $L^2$-orthogonal.
\end{prop}
In light of Theorem~\ref{harmonic}(a), $\gamma$ is equivalently an $L^2$
harmonic form.

When $(M,\Phi)$ is an asymptotically cylindrical $\Spin(7)$-manifold
the Hodge Laplacian maps the sections $\Omega^r_k(M)$ of each
$\Spin(7)$-invariant subbundle $\Lambda^r_k$ of $r$-forms~\eqref{rep4} into
itself \cite[\S 3.5]{joyce}. Then there is a well-defined decomposition
$\HH^r(M)=\oplus_k\HH^r_k(M)$ into spaces of bounded harmonic forms
$\HH^r_k(M)\subset\Omega^r_k(M)$.

The Levi--Civita connection of $g(\Phi)$ on $\Lambda^r T^* M$ also
preserves the decomposition~\eqref{rep4} determined by a parallel
$\Spin(7)$-structure $\Phi$. Let $\widehat\HH^r_k(M)$ denote the space of
parallel forms in $\Omega^r_k(M)$. The following result will be
useful in the next section.

\begin{prop}\label{2-forms}
Let $M$ be an asymptotically cylindrical $\Spin(7)$-manifold with
cross-section $Y$ and with the $\Spin(7)$-structure $\Phi$ asymptotic to
$\Phi_\infty=dt\we\phi+*_\phi(\phi)$ for a torsion-free $G_2$-structure
$\phi$ on~$Y$.
Suppose that $b^1(M)=0$. Then $\widehat{\HH}^2_7(M)=\HH^2_7(M)$ and the map
\begin{equation}\label{atinfty}
\alpha\in \widehat{\HH}^2_7(M) \mapsto \alpha_\nu\in\HH^1(Y)
\end{equation}
is linear isomorphism, where $\alpha_\nu$ is defined by~\eqref{asympt}.
\end{prop}
\begin{pf}
Every $\alpha\in\Omega^2_7(M)$ on a $\Spin(7)$-manifold $M$ satisfies
the Weitzenb\"ock formula
\begin{equation}\label{W}
\Delta\alpha=\nabla^*\nabla\alpha.
\end{equation}
(\cite[Prop. 10.6.5]{joyce}), so every
$\alpha\in\widehat{\HH}^2_7(M)$ is a bounded harmonic form,
$\alpha\in\HH^2_7(M)$.

Conversely, let $\alpha\in\HH^2_7(M)$. Then by Theorem~\ref{harmonic}
$\alpha$ is asymptotic to a translation-invariant harmonic form
$\alpha_\infty=dt\we\alpha_\nu+\alpha_\tau$ on the cylindrical end
$M_\infty$. Let
$\Lambda^2T^*M_\infty=(\Lambda^2_7)_\infty\oplus(\Lambda^2_{21})_\infty$
denote the $\Spin(7)$-invariant `type decomposition' of 2-forms
induced by~$\Phi_\infty$. As $\alpha\in\Omega^2_7(M)$ with respect to
$\Phi$ and $\Phi$ is asymptotic to $\Phi_\infty$, the
$(\Lambda^2_{21})_\infty$ component of $\alpha$ is decaying to zero as
$t\to\infty$. On the other hand, both $\alpha_\infty$ and $\Phi_\infty$
are translation-invariant and so is the $(\Lambda^2_{21})_\infty$ component
of~$\alpha_\infty$. Therefore, $\alpha_\infty$ is a section of
$(\Lambda^2_7)_\infty$.

In terms of the $G_2$-structure $\phi$ and the forms on~$Y$, the bundles
$(\Lambda^2_j)_\infty$ can be written as
\begin{gather}
(\Lambda^2_7)_\infty =
\{dt\we *_\phi((*_\phi\phi)\we v\lrcorner\phi)+
3v\lrcorner\phi \;|\; v\in TY\},
\label{type.7}
\\
(\Lambda^2_{21})_\infty =
\{dt\we *_\phi((*_\phi\phi)\we\alpha)-\alpha \;|\;
\alpha\in\Lambda^2 T^* Y\}
\label{type.21}
\end{gather}
(cf. \cite[Theorem~11.4]{salamon-walpuski}), where we also noted that a
$G_2$-structure $\phi$ induces a $G_2$-invariant rank 7 subbundle of
2-forms $\Lambda^2_7T^* Y=\{v\lrcorner\phi : v\in TY\}$. The map 
\begin{equation}\label{iso7}
v\lrcorner\phi  \mapsto *_\phi(*_\phi(\phi)\we v\lrcorner\phi)
\end{equation}
defines a $G_2$-equivariant isometry of bundles $\Lambda^2_7T^* Y\to T^*Y$
and gives an isomorphism between subspaces of parallel forms in
$\Lambda^2_7T^* Y$ and~$T^*Y$. As the $G_2$-manifold $Y$ is Ricci-flat the
1-forms on~$Y$ satisfy the Weitzenb\"ock formula~\eqref{W}, without
curvature terms. Then the integration by parts on a compact $Y$ shows
that $\HH^1(Y)$ is precisely the space of parallel 1-forms.

We obtain that the leading asymptotic term $\alpha_\infty$ of $\alpha$
is linearly determined by $\alpha_\nu$ via~\eqref{type.7} and
$\alpha_\infty$ is parallel as $\alpha_\nu\in\HH^1(Y)$. Then $\nabla\alpha$
decays to zero along the end of~$M$, so the integration by parts argument
is valid for $\alpha$ and shows that $\nabla\alpha=0$. It follows that
$\HH^2_7(M)=\widehat\HH^2_7(M)$ and \eqref{atinfty} is a well-defined
injective linear map.

For the surjectivity, it suffices to show that the dimension of the
space of parallel forms in $\Omega^2_7(M)$ is $b^1(Y)$.
We have $\dim\HH^2(M)=b^1(Y)+b^2(M)$ as the hypothesis $b^1(M)=0$ yields
$b^2_c(M)=b^1(Y)+b^2_0(M)$ from the exact sequence~\eqref{exact}.
On the other hand, $\dim\HH^2(M)=\dim\widehat\HH^2_7(M)+\dim\HH^2_{21}(M)$
from previous arguments. We see from~\eqref{type.21} and
Theorem~\ref{harmonic} that the asymptotic term $\alpha_\infty$ for each
$\alpha\in\HH^2_{21}(M)$ is determined by $\alpha_\tau\in\HH^2(Y)$ and if
$\alpha_\tau=0$ then $\alpha$ is an $L^2$ form. It follows find that
$\HH^2_{21}(M)$ can be identified with the direct sum of the images of
homomorphisms $H^2(M)\to H^2(Y)$ and $H^2_c(M)\to H^2(M)$ in~\eqref{exact}
and therefore $\dim\HH^2_{21}(M)\le b^2(M)$. As $\dim\widehat\HH^2_7(M)\le
b^1(Y)$ by the above work, we find that the equalities are attained,
$\dim\HH^2_{21}(M)=b^2(M)$ and $\dim\widehat\HH^2_7(M)=b^1(Y)$.
\end{pf}

The next result would be needed in \S\ref{de.Rham} below.

\begin{prop}\label{3-forms}
Let $M$ be an asymptotically cylindrical $\Spin(7)$-manifold with
cross-section $Y$ and with the $\Spin(7)$-structure $\Phi$ asymptotic to
$\Phi_\infty=dt\we\phi+*_\phi\phi$ for a torsion-free $G_2$-structure
$\phi$ on~$Y$. Suppose that $b^1(M)=0$.

Then $\HH^3_8(M)=\{0\}$ and for each $\alpha\in\HH^3(M)$ the asymptotic
component $\alpha_\tau\in\HH^3(Y)$ defined by~\eqref{asympt} is orthogonal
to $\phi$ at each point in~$Y$ with the metric $g(\phi)$.
\end{prop}
\begin{pf}
The vector bundles $T^*M$ and $\Lambda^3_8T^*M$ are associated to
isomorphic representations of $\Spin(7)$ and it follows from Weitzenb\"ock
formula that $\HH^3_8(M)\cong\HH^1(M)$.
On the other hand, as $b^1(M)=0$ it follows from Corollary~\ref{bounded}
and the exact sequence~\eqref{exact} that $\HH^1(M)=\{0\}$.
Thus $\HH^3_8(M)=\{0\}$ and $\HH^3(M)=\HH^3_{48}(M)$.

For the second claim, an argument similar to that in
Proposition~\ref{2-forms} shows that for each $\alpha\in\HH^3_{48}(M)$
the asymptotic term $\alpha_\infty$ is in $\HH^3_{48}(\RE\times Y)$
relative to the cylindrical $\Spin(7)$-structure $\Phi_\infty$.
That $\alpha_\tau$ is pointwise orthogonal to $\phi$ on~$Y$ now follows
directly from the determination of $\Spin(7)$-invariant subbundle
$\Lambda^3_{48}T^*(\RE\times Y)$ in terms of the $G_2$-invariant subbundles
of $\Lambda^r_kT^*Y$ given in~\cite[Theorem~8.4 and 11.4]{salamon-walpuski}.
\end{pf}

\section{Holonomy and deformations of asymptotically
cylindrical $\Spin(7)$-manifolds} 
\label{proper.spin7}

The holonomy group of $\Spin(7)$-manifold can be a proper
subgroup of $\Spin(7)$, even when the manifold is compact and
simply-connected; one source of examples is provided by Calabi--Yau
manifolds. On the other hand, the following general criterion was proved
by Bryant.

\begin{thrm}[{\cite[cf.\ pages~565--566]{bryant}}]\label{spin7crit}
Let $M$ be a simply-connected $\Spin(7)$-manifold. Then the holonomy of $M$
is $\Spin(7)$ if and only if $M$ admits no non-zero parallel 1-forms or
2-forms.
\end{thrm}

For the asymptotically cylindrical $\Spin(7)$-manifolds, the holonomy is
determined topologically from Betti numbers of the cross-section as the
next result shows.

\begin{thrm}\label{hol8d}\label{irred.cor}
If $M$ is a simply-connected asymptotically cylindrical
$\Spin(7)$-manifold with non-empty cross-section~$Y$ then the holonomy
group of~$M$ is one of $\Spin(7)$, $G_2$, $SU(4)$ and $SU(2)\times SU(2)$.

(a)
The holonomy of~$M$ is $\Spin(7)$ if and only if
$b^0(Y)=1$ (i.e. $M$ has a single end) and $b^1(Y)=0$.

(b)
The holonomy of~$M$ is $G_2$ if and only if $M$ is a cylinder $Y\times\RE$
with a product metric, where $Y$ is a compact simply-connected
7-manifold with holonomy $G_2$.

(c) The holonomy of~$M$ is $SU(4)$ if and only if $b^1(Y)=1$.

(d) The holonomy of~$M$ is $SU(2)\times SU(2)$ if and only if $b^1(Y)=3$.
In this case, $M=S_1\times S_2$ is the Riemannian product of a
Ricci-flat K\"ahler K3 surface $S_1$ and a simply-connected asymptotically
cylindrical Ricci-flat K\"ahler surface $S_2$ with cross-section
homeomorphic to the 3-torus.
\end{thrm}
\begin{rmk}
We construct examples realizing the possibility (a) in this paper.
The case (b) is implemented by any simply-connected $G_2$-manifold,
e.g.\ \cite{joyce} or \cite{g2paper}. Regarding (c) and (d), a method
producing examples of asymptotically cylindrical Calabi--Yau metrics is
given in \cite{g2paper}. See also \cite{HHN} and \cite{biquard-minerbe}.
It is shown in \cite{hein} that the flat 3-torus cross-section at infinity
of $S_2$ in (d) need not in general be isometric to a Riemannian product of
$S^1$ and a 2-torus.
\end{rmk}

Much of Theorem~\ref{hol8d} is known in some form. The criterion (a) is
a direct consequence of the result of
Nordstr\"om~\cite[Theorem~4.1.19]{jn-thesis} proved using spinors and
(c) is obtained, again by spinor methods, in~\cite[Theorem~B]{HHN} which
also excludes the holonomy $Sp(n)$.

We give a self-contained argument for a simply-connected~$M$ using
the octonions algebra and parallel differential forms without appealing to
spinors or holomorphic geometry.

\begin{pf}[Proof of Theorem~\ref{hol8d}]
As every $\Spin(7)$-manifold is Ricci-flat, the Cheeger--Gromoll splitting
theorem \cite{split.thm} implies that a connected asymptotically cylindrical
$\Spin(7)$-manifold $M$ has either one end or two ends (see also~\cite{salur}).

Furthermore, if $M$ has two ends, then it is necessarily a cylinder
$M=\RE\times Y$ with a product metric. Then $Y$ is a compact
simply-connected $G_2$-manifold and by~\cite[Prop. 10.2.2]{joyce}
the holonomy group of $Y$, hence also of~$M$, is~$G_2$.
Conversely, if the holonomy of $M$ is~$G_2$, then an application of the
de Rham theorem~\cite{deRham} shows that $M$ is the Riemannian product of a
complete $G_2$-manifold $Y$ and a 1-manifold. As $M$ is simply-connected
with (only) cylindrical ends, we deduce that $M=\RE\times Y$ and $Y$ is
compact simply-connected with holonomy~$G_2$. This completes the proof of~(b).

Suppose now that $M$ has a single end. Applying the de Rham theorem again
we find that the holonomy group of~$M$ cannot act trivially on any subspace
of~$\RE^8$ of positive dimension. The only connected subgroups of
$\Spin(7)$ with the latter property and allowed by Berger's classification
are $\Spin(7)$, $SU(4)$, $Sp(2)$ and $SU(2)\times SU(2)$.

Let $b^1(Y)=0$. If the holonomy of $M$ is not $\Spin(7)$,
then it is a subgroup of $SU(4)$ and $M$ has a parallel 2-form, the
K\"ahler form $\omega$. But as we saw in \eqref{kaehler} then $\omega$
is in $\Omega^2_7(M)$ and we obtain $b^1(Y)>0$ by Proposition~\ref{2-forms},
a contradiction. Conversely, if the holonomy of $M$ is $\Spin(7)$, then by
Theorem~\ref{spin7crit} there are no non-zero parallel 2-forms, hence
$b^1(Y)=0$ by Proposition~\ref{2-forms}. This completes the proof of~(a).

If the holonomy of~$M$ is contained in $SU(4)$, then $M$ has a parallel
complex structure~$I$ and there is an $SU(4)$-invariant decomposition of
real 2-forms on $M$ refining the $\Spin(7)$-invariant
decomposition~\eqref{rep4},
\begin{equation}\label{su4spin7forms}
\Lambda^2_7=\Lambda^2_1\oplus\Lambda^2_6,
\qquad
\Lambda^2_{21}=\Lambda^2_{15}\oplus\Lambda^2_6,
\end{equation}
according to the types (1,1) and $(2,0)\oplus(0,2)$ with respect to~$I$.
More explicitly, the fibres of the two $\Lambda^2_6$ subbundles
in~\eqref{su4spin7forms} are, respectively, the $\pm 2$-eigenspaces of
$*(\cdot\we\Re\theta_0)$, where $\theta_0$ is the $(4,0)$ form
in~\eqref{su4spin7}, cf. \cite[Theorem~11.4]{salamon-walpuski}.
These two rank 6 subbundles are isomorphic, being associated to isomorphic
representations of $SU(4)$. The sections of $\Lambda^2_6$ satisfy the
Weitzenb\"ock formula~\eqref{W} with vanishing curvature term and the
argument of Proposition~\ref{2-forms} applies to show that the bounded
harmonic sections of each $\Lambda^2_6$ are parallel.

If the holonomy of $M$ is exactly $SU(4)$ then the only parallel 2-forms
on~$M$ are constant sections of $\Lambda^2_1$ (constant multiples of the
K\"ahler form) and $b^1(Y)=1$.

If the holonomy of $M$ is a proper subgroup of $SU(4)$, then it is contained
in $Sp(2)$ and $M$ admits three parallel  complex structures and three
respective K\"ahler forms which are sections of $\Lambda^2_7$. These
K\"ahler forms are point-wise orthogonal to each other, thus two of 
them are sections of $\Lambda^2_6\subset\Lambda^2_7$ and, respectively,
there are further two parallel sections of $\Lambda^2_6\subset\Lambda^2_{21}$.
We thus obtain a 5-dimensional space of parallel 2-forms on~$M$.
As the holonomy representation of $Sp(2)$ stabilizes only three linearly
independent 2-forms on $\RE^8$, we deduce that the holonomy of $M$ reduces
to a proper subgroup of $Sp(2)$ and therefore must be $SU(2)\times SU(2)$.

Now if the holonomy of $M$ is $SU(2)\times SU(2)$, then applying the de
Rham theorem again we obtain that $M=S_1\times S_2$ is a Riemannian product
of two 4-manifolds with holonomy $SU(2)$, i.e.\ $S_i$ are Ricci-flat K\"ahler
surfaces with a trivial canonical bundle. Both $S_i$ are simply-connected
and one of these, $S_1$ say, is compact, hence a K3 surface, and the metric
on $S_2$ is asymptotically cylindrical. The 3-dimensional cross-section of
$S_2$ at infinity then is a flat 3-torus, thus $b^1(Y)=3$ in this case.
This completes the proof of Theorem~\ref{hol8d}.
\end{pf}

Joyce investigated the deformation theory of torsion-free
$\Spin(7)$-structures on a compact 8-manifold $M$ and proved that  these
have a smooth moduli space of dimension $b^4_-(M)+b^2_7(M)+1$, which
becomes $b^4_-(M)+1$ if $M$ has full holonomy $\Spin(7)$. (Here $b^r_j(M)$
denote the dimensions of the spaces of harmonic forms of respective types.)
Cf.~\cite[Theorem~10.7.1]{joyce}. 

Nordstr\"om obtained an analogue of Joyce's result for 8-manifolds with
cylindrical ends, where the moduli space is defined using diffeomorphisms
isotopic to the identity and asymptotically to those preserving the 
product cylindrical structure of the ends. Recall from \S\ref{as.cyl}
the definition of subspaces $H^r_0$ in the de Rham cohomology.

\begin{thrm}[{\cite[Theorem 4.3.2 and Prop. 4.3.3]{jn-thesis}}]\label{moduli}
Let $(M,\Phi)$ be a connected asymptotically cylindrical
$\Spin(7)$-manifold. Then the moduli space of torsion-free asymptotically
cylindrical $\Spin(7)$-structures on~$M$ is a smooth manifold
of dimension $b^4(M)-b^4_+(M)-b^1(M)+b^1(Y)+1$, where $b^4_+(M)$ is the
dimension of a maximal positive subspace for the cup-product on $H^4_0(M)$.
\end{thrm}

For a simply-connected asymptotically cylindrical~$M$ with holonomy equal
to $\Spin(7)$, the dimension of the moduli space therefore is
$b^4(M)-b^4_+(M)+1$. It can be checked using the argument of
Corollary~\ref{bounded} and the Poincar\'e duality that $b^4(M)-b^4_+(M)$
is the dimension of the space $\HH^4_-(M)\subset\Omega^4_-(M)$ of bounded
harmonic anti-self-dual 4-forms on~$M$.

\section{Existence of asymptotically cylindrical torsion-free
  $\Spin(7)$-structures}
\label{exist}

In order to construct asymptotically cylindrical metrics with holonomy
$\Spin(7)$ on 8-manifolds with cylindrical ends, we shall obtain
$\Spin(7)$-structures (more precisely, 1-parameter families) with arbitrary
small torsion, satisfying Theorem~\ref{exist.spin7} which is the main
result of this section. We modify Joyce's existence result
\cite[Theorem~13.6.1]{joyce} (and Proposition 15.2.13 {\it op.cit.}) for
torsion-free $\Spin(7)$-structures on compact manifolds, extending it to
manifolds with cylindrical ends. The strategy of the proof of
Theorem~\ref{exist.spin7} has some analogy with \cite[Theorem~3.1]{KN} for
asymptotically cylindrical torsion-free $G_2$-structures on 7-manifolds.

\begin{thrm}\label{exist.spin7}
Let $\lambda$, $\mu$, $\nu$ be positive constants.
Then there exist positive constants $\kappa$, $K$ and $K'$ such that
whenever $0<s\le\kappa$ the following is true.

Let $M$ be an 8-manifold with cylindrical end~$M_\infty=\RE_+\times Y$.
and suppose that an admissible 4-form $\tPhi$ defines on $M$ a
$\Spin(7)$-structure which is asymptotically cylindrical with
rate~$\delta_M$.
Suppose that $\tPhi$ is closed on $M_\infty$ and there is a smooth compactly
supported 4-form $\psi$ on~$M$ satisfying $d\tPhi+d\psi=0$. If
\begin{enumerate}
\itemsep0em
\item
\begin{equation}\label{psi}
\|\psi\|_{L^2}\le\lambda s^{13/3}
\text{ and }
\|d\psi\|_{L^{10}}\le\lambda s^{7/5},
\end{equation}
\item\label{inj}
the injectivity radius $\delta(\tPhi)$ of the metric $g(\tPhi)$
satisfies $\delta(\tPhi)\ge\mu s$, and
\item\label{curv}
the Riemann curvature $R(\tPhi)$ of the metric $g(\tPhi)$ satisfies
$\|R(\tPhi)\|_{C^0}\le\nu s^{-2}$.
\end{enumerate}
then there is a smooth anti-self-dual 4-form $\eta$ on $M$, exponentially
decaying to zero with all derivatives as $t\to\infty$, such that 
$\|\eta\|_{C^0}\le Ks^{1/3}$ and
\begin{equation}\label{spin7.eqn}
d\eta=d\psi+dF(\eta)
\end{equation}
and $\Phi=\Theta(\tPhi + \eta)$ is a torsion-free $\Spin(7)$-structure,
with $\|\Phi-\tPhi\|_{C^0}\le K's^{1/3}$.
The exponential decay rate of $\eta$ can be taken to be
$\delta=\frac12\min\{\delta_M,\delta_Y\}$, where $\delta_Y^2$ is the
smallest positive eigenvalue of the Hodge Laplacian on $\Omega^*(Y)$ for
the metric of the cross-section of~$M$ at infinity.
\end{thrm}

The difference between Theorem~\ref{exist.spin7} and
\cite[Theorem~13.6.1]{joyce} is that in the present case $M$ is non-compact
with a cylindrical end and we added appropriate hypotheses on the
asymptotical behaviour of $\tPhi$ and $\psi$ on the end and are claiming
an asymptotically cylindrical property of the torsion-free $\Phi$.
Formally, Theorem~\ref{exist.spin7} includes the statement of
\cite[Theorem~13.6.1]{joyce}, corresponding to the case when $Y=\emptyset$
(so $M$ is compact). 

In the rest of this section we prove Theorem~\ref{exist.spin7}. Firstly, we
note that it suffices to show the existence of exponentially decaying
anti-self-dual smooth 4-form $\eta$, with `small' uniform norm,
satisfying~\eqref{spin7.eqn}. The argument is local and readily carries
over from \cite[Theorem~13.6.1]{joyce}. We then adapt the contraction
mapping method used by Joyce \cite[\S 13.7]{joyce} on compact 8-manifolds
to find an anti-self-dual form $\eta$ solving~\eqref{spin7.eqn} on
asymptotically cylindrical $M$ and use elliptic regularity to show that 
the solution $\eta$ is smooth. In fact, the method will also show that
$\eta$ decays along the cylindrical end $\RE_+\times Y\subset M$ at
$t\to\infty$, uniformly in the $Y$ variables.
Our last claim that the solution decays at an {\em exponential} rate
has no analogue in the compact case. To prove this claim we interpret 
\eqref{spin7.eqn} on the end of~$M$ as an infinite-dimensional flow in~$t$
and use results about the hyperbolic stationary points.

Suppose that $\eta\in\Omega^4_-(M)$ satisfies~\eqref{spin7.eqn} and 
$\|\eta\|_{C^0}\le Ks^{1/3}$ with $s$ sufficiently small. Then, by the
discussion in~\S\ref{synopsis}, $\Phi=\Theta(\tPhi + \eta)$ is a
well-defined $\Spin(7)$-structure and $\Phi=\tPhi + \eta - F(\eta)$ as
$\eta\in\Omega^4_-=\Omega^4_{35}$. Therefore,
$d\Phi=-d\psi+d\eta-dF(\eta)=0$ and $\Phi$ is torsion-free. Also,
$\|\Phi-\tPhi\|_{C^0}\le|\eta|+|F(\eta)|<K' s^{1/3}$ noting the
estimate~\eqref{quadratic} for $F$.

Thus we wish to find an appropriate solution of~\eqref{spin7.eqn}.
The following is a version of the contraction mapping result
\cite[Prop. 13.7.1]{joyce} adapted to the asymptotically cylindrical
base manifold.
\begin{prop}\label{cmp}
Assume that an 8-manifold $M$ with cylindrical end and a $\Spin(7)$
structure $\Phi$ on~$M$ satisfy all the hypotheses of
Theorem~\ref{exist.spin7}.

Then there exist positive constants $\kappa$, $C$ and $K$ depending only
on $\lambda$ in~\eqref{psi} such that if $s\le\kappa$ then there 
exists a sequence $\eta_j\in \Omega^4_-(M)$, $\eta_0=0$, of smooth forms
exponentially decaying with all derivatives with rate $\delta$ (defined in
Theorem~\ref{exist.spin7}) and satisfying, for each $j=1,2,\ldots$, the
equation
\begin{equation}\label{iterate}
d\eta_j=d\psi+dF(\eta_{j-1})
\end{equation}
and the inequalities
\begin{enumerate}
\itemsep0em
\item\label{L2}
$\|\eta_j\|_{L^2}\le 4\lambda s^{13/3}$,
\item
$\|\nabla\eta_j\|_{L^{10}}\le C s^{2/15}$,
\item\label{unif}
$\|\eta_j\|_{C^0}\le Ks^{1/3}$
\item\label{L2d}
$\|\eta_j-\eta_{j-1}\|_{L^2}\le 4\lambda 2^{-j} s^{13/3}$,
\item
$\|\nabla(\eta_j-\eta_{j-1})\|_{L^{10}}\le C 2^{-j} s^{2/15}$,
\item
$\|\eta_j-\eta_{j-1}\|_{C^0}\le K 2^{-j} s^{1/3}$.
\end{enumerate}
\end{prop}
\begin{pf}
The proof proceeds by induction on~$j$. Assume that there exist
$\eta_0,\ldots,\eta_{j-1}$ satisfying all the assertions of
Proposition~\ref{cmp}. Then $\eta_{j-1}\in\Omega^4_-(M)$ is exponentially
decaying and satisfies the uniform estimate \ref{unif} with small~$s$, so
$F(\eta_{j-1})$ is well-defined and satisfies the quadratic
estimate~\eqref{quadratic}. Therefore, $\psi+F(\eta_{j-1})$ decays with the 
same exponential rate as $\eta_{j-1}$. Applying Proposition~\ref{decomp} we
can write $\psi+F(\eta_{j-1})=d\alpha+d^*\beta+\gamma$ and
$$
d\psi+dF(\eta_{j-1})=dd^*\beta=d(d^*\beta-*d^*\beta+\tgamma),
$$
for some exponentially decaying $d^*\beta$ and an arbitrary $L^2$ harmonic
form $\tgamma$. By Theorem~\ref{harmonic}(a)), $\tgamma$ is closed and may
be chosen to represent any given cohomology class in $H^4_0(M)$. We choose
$\tgamma$ so that $\tgamma-\gamma)\in\HH^4_+(M)$, where
we denoted by $\HH^4_\pm=L^2\Omega^4_\pm(M)\cap\HH^4(M)$ the spaces of
self- and anti-self-dual $L^2$ harmonic 4-forms. Note that the forms in
$\HH^4_\pm$ represent exactly the classes in a maximal positive,
respectively negative, subspace $H^4_\pm(M)$ for the cup-product on
$H^4_0(M)$. Now define $\eta_j:=d^*\beta-*d^*\beta+\gamma$ and this is a
smooth anti-self-dual 4-form exponentially decaying with all derivatives
and satisfying~\eqref{iterate}. 

The method used in \cite[\S 13.7]{joyce} to prove the estimates (a)--(f)
for a compact base manifold $M$ applies {\it mutatis-mutandis} to the
present situation. The role of harmonic 4-forms and their cohomology
classes in the compact case is now taken by the exponentially decaying
(equivalently, $L^2$) harmonic 4-forms on~$M$ and the subspace
$H^4_0(M)\subset H^4(M)$. The respective integrals in the proof of \ref{L2}
and \ref{L2d} are finite on~$M$ as each $\eta_j$ decays exponentially with
all derivatives along the cylindrical end. The proof of (b),(c),(e),(f) in
\cite[\S 13.7]{joyce} uses estimates for Sobolev norms on complete
Riemannian 8-manifolds (it is at this point that the hypotheses~\ref{inj}
and \ref{curv} of Theorem~\ref{exist.spin7} are required) and the quadratic
estimates \eqref{quadratic},\eqref{quadratic2} for~$F$ and carries over
without changes. We therefore omit the details.
\end{pf}

When $s$ is sufficiently small the sequence $\{\eta_j\}$ given by
Proposition~\ref{cmp} is Cauchy and converges in the $L^{10}_1$ norm and
also in the $L^2$ and $C^0$ norms on~$M$, to $\eta\in\Omega^4_-(M)$ say.
The limit $\eta$ is a weak $L^{10}_1$ solution to the~\eqref{spin7.eqn} and
satisfies
$$
\|\psi\|_{L^2}\le K s^{13/3},
\qquad
\|\psi\|_{C^0}\le K s^{13/3}
\quad
\|\nabla\psi\|_{L^{10}}\le K s^{7/5},
$$
for some $K>0$ independent of $0<s<\kappa$.

\begin{prop}\label{regularity}
If $s>0$ is sufficiently small, then

(i) $\eta$ is in $L^{10}_k(M)$, for each $k>0$, therefore smooth on (M), and

(ii) $\eta$ on the cylindrical end of~$M$ converges to zero with all
derivatives, uniformly on~$Y$, as $t\to\infty$. 
\end{prop}
\begin{pf}
For part (i), we note that the exterior derivative on $\Omega^4_-(M)$ is an
overdetermined elliptic operator, i.e.\ has injective symbol, and the
standard interior elliptic $L^p$ estimates still hold for it. (More
explicitly, it is a restriction of the elliptic operator arising from
elliptic complex $0\to\Omega^4_-(M)\to\Omega^5(M)\to\ldots$.)
As $F(\eta)$ satisfies a quadratic estimate we can write the term
$dF(\eta)$ in~\eqref{spin7.eqn} as $G(\eta,\nabla\eta)$, where $G$ is 
a smooth function, linear in the second argument and $G(0,\cdot)\equiv 0$.
Then equation satisfied by $\eta$ takes the form
$$
d\eta-G(\eta,\nabla\eta)=d\psi.
$$
We regard the left-hand side as a linear operator in $P_\zeta(\eta):=
d\eta-G(\zeta,\nabla\eta)$ with $\zeta=\eta$. Assuming the $\kappa$ is
sufficiently small, we have the $C^0$ norm of $\zeta$ sufficiently small.
Then $P_\zeta$ is overdetermined elliptic as this is an open condition.

Is $\eta\in L^{10}_k(M)$ for some $k\ge 0$, then the coefficients of
$P_\zeta$ are also in $L^{10}_k(M)$, hence in $C^{k,1/5}(M)$ by Sobolev
embedding. Applying \cite[Theorems 6.2.5, 6.2.6]{morrey} we deduce that
$\eta$ is locally in $L^{10}_{k+1}(M)$. Furthermore, the estimate of
the $L^{10}_{k+1}$ norm on local neighbourhoods can be taken uniformly on
$M$ as the metric on~$M$ is asymptotically cylindrical. It follows that 
$\eta\in L^{10}_{k+1}(M)$ and then in $C^\infty(M)$ by induction on~$k$.

Part (ii) follows directly from the argument of \cite[Cor. 3.7]{KN} which 
exploits the finiteness of $L^{10}_k$ norm of $\eta$ and a uniform bound
for the local embedding $L^{10}_{k+1}\subset C^k$ on geodesic balls of fixed
radius in~$M$.
\end{pf}

The exponential rate of decay requires additional work. We 
consider the `Spin(7) equation' \eqref{spin7.eqn} on the cylindrical end
$\RE_+\times Y\subset M$ as an evolution equation on~$Y$ with the zero
solution corresponding to a stationary point. The decaying solution then
may be considered as a path of forms parameterised by $t$ on the
cross-section 7-manifold $Y$. The key property is that the stationary point
zero is {\em hyperbolic}.  The exponential decay property is well-known for
hyperbolic stationary points of flows in finite-dimensional spaces. For an
infinite-dimensional setting, a suitable argument is given in \cite[Chap.~5,
Theorem~5.2.2]{MMR}. In fact, the result is slightly more general and
asserts that every solution path converging to a stationary point is
exponentially asymptotic to a path in an invariant `centre manifold' ---
when the stationary point is hyperbolic, a centre manifold is just this
stationary point.

For each $p$, the space of $p$-forms on~$\RE_+\times Y$
is naturally isomorphic to the one-parameter families
$\RE\to(\Omega^{p-1}\oplus\Omega^p)(Y)$. For anti-self-dual 4-forms
on~$\RE_+\times Y$ we may write
\begin{equation}\label{eta.cyl}
\eta|_{M_\infty}= dt\we(\eta_3(t)+\wp(t,x,\eta_3)- *_Y\eta_3(t),
\end{equation}
where $\wp(t,x,\eta_3)$ is linear in $\eta_3$, smooth in $x\in Y$ and
$O(e^{-\delta_M t)}$ with all derivatives, for $t>0$.

As $\psi$ is compactly supported, we may take $\psi=0$ on $M_\infty$, so
the equation satisfied by $\eta$ becomes
\begin{equation}\label{on.cyl}
d\eta=dF(\eta).
\end{equation}
We know that $\eta_3(t)\to 0$ uniformly on~$Y$ as $t\to\infty$.

Write $F(\eta)=dt\we F_3(\eta)-F_4(\eta)$ and then \eqref{quadratic}
implies that $F_3,F_4$ on $\{t\}\times Y$ satisfy quadratic estimates in
$\eta_3(t)$ for each~$t$,
$$
|F_i(\eta'_3)-F_i(\eta''_3)|<
C_2|\eta'_3-\eta''_3|(|\eta'_3|+|\eta''_3|),\qquad  i=3,4,
$$
with $C_2$ independent of~$t$.

Differentiating~\eqref{eta.cyl} and $F(\eta)$, plugging in~\eqref{on.cyl}
$$
-dt\we(d_Y\eta_3+\frac\p{\p t}(*_Y\eta_3))-d_Y*_Y\eta_3=
-dt\we(d_Y F_3+\frac\p{\p t} F_4)-d_Y F_4
$$
and equating the components, we obtain
\begin{equation}\label{flow}
d_Y(*_Y\eta_3-F_4)=0,\qquad
\frac\p{\p t}(*_Y\eta_3-F_4)=d_Y (F_3-\eta_3-\wp),
\end{equation}
so $*_Y\eta_3(t)-F_4(\eta_3(t))$ is exact on $\{t\}\times Y$ for each $t>0$.

The function assigning to $\eta_3\in\Omega^3(Y)$ the co-closed component of
$*_Y\eta_3-F_4(\eta_3)$, in the Hodge decomposition on compact~$Y$, is smooth 
and has surjective derivative at $\eta_3=0$. By the implicit function
theorem in Banach spaces the set
$$
\mathfrak{Y}_\eps=\{\eta_3\in L^2\Omega^3(Y):\|\eta_3\|<\eps,\;
*_Y\eta_3-F_4(\eta_3) \text{ is exact}\}
$$
for some $\eps>0$ is a graph of smooth function
$d^*_Y \bigl(L^2_1\Omega^4(Y)\bigr)\to \Ker d\cap L^2\Omega^3(Y)$.
In particular $\mathfrak{Y}_\eps$ is a Banach submanifold of
$L^2\Omega^3(Y)$ with the tangent space at $\eta_3=0$
$$
T_0\mathfrak{Y}_\eps=d^*_Y \bigl(L^2_1\Omega^4(Y)\bigr).
$$
By the Hodge theory, the space of co-exact forms
$T_0\mathfrak{Y}_\eps$ for a compact $Y$ is a closed subspace of $L^2$
3-forms, thus a Banach space. Further, the $L^2$-orthogonal projection
$\pr_Z:L^2\Omega^3(Y)\to T_0\mathfrak{Y}_\eps$ is a bounded
linear map. Denote $\tilde\eta_3:=\pr_Z(\eta_3)$.
The restriction of $\pr_Z$ to $\mathfrak{Y}_\eps$ is invertible with
$\eta_3=\tilde\eta_3-\tilde{F}_3(\tilde\eta_3)$ and $\tilde{F}_3$
satisfying a quadratic estimate.

Substituting the latter formula in the second equation of~\eqref{flow}
we obtain that $\tilde\eta_3$ satisfies a flow equation in on a
neighbourhood of 0 in the space of co-exact $L^2$ 3-forms on~$Y$
\begin{equation}\label{exact.flow}
\frac{\p\te}{\p t}=
-*_Yd_Y\te + Q(\te) + \tilde\wp(t,\te),
\end{equation}
with $Q(\tilde\eta_3)$ of quadratic order in $\tilde\eta_3$ and
$\tilde\wp$ of linear order in $\tilde\eta_3$
and $\tilde\wp(t,\te)=O(e^{-\delta_M t)}$ for $t>0$.

The operator $*_Y\,d_Y$ defines a linear isomorphism of the space
$d^*_Y\Omega^4(Y)$ onto itself and is formally self-adjoint. Its square is
the (restriction of) Hodge Laplacian. The space of co-exact 3-forms splits
into $L^2$-orthogonal direct sum $d^*_Y \Omega^4(Y)=T_+\oplus T_-$ of
positive and negative eigenspaces of $-*_Y\,d_Y$ and $\eta_3=0$ is thus
a hyperbolic stationary point of the infinite-dimensional
flow~\eqref{exact.flow}.

\begin{prop}\label{exp.decay}
Let $\te(t)$, $t\ge t_0$, be a path of co-exact 3-forms on~$Y$
satisfying the equation~\eqref{exact.flow} and such that $\tilde\eta_3(t)\to 0$
uniformly on~$Y$ as $t\to\infty$.

Then $\|\tilde\eta_3\|_{L^p_k(\{t\}\times Y)}<C_{p,k}e^{-(\delta/2) t})$ for
$t\ge t_0$ (here $p>1$, $k=0,1,2,\ldots$) with $\delta$ defined in
Theorem~\ref{exist.spin7}.
\end{prop}
\begin{pf}
This is an instance of a general property for partially hyperbolic fixed
point of flows. Our argument is similar to that in
\cite[\S 5.4]{MMR}, but in the present case the linear part
has no pure imaginary eigenvalues and some details are simplified.
It will suffice to consider $p=2$ as any given $L^p_j$ norm of a smooth
$\te$ on~$Y$ is controlled by its $L^2_k$ norm for a sufficiently large~$k$.

The hypothesis implies that there is, for
each $E_0>0$ and $k=1,2,\ldots$, a $T_{0,k}>0$, so that 
$\sup_{t\in [T_{0,k},\infty)}\|\te\|_{L^2_k(\{t\}\times Y)}<E_0$.
Given a $\te\in T_0\mathfrak{Y}_\eps$ we can write $\te=\te^+ + \te^-$, for
some unique $\te^+\in T^+$ and $\te^-\in T^-$.

We claim that for a sufficiently small $E_0$, the paths $\te^\pm(t)$
satisfy the differential inequalities
\begin{subequations}\label{claim.ineq}
\begin{align}
\frac{\p}{\p t}\|\te^+\|_{L^2_k}&\ge
\delta(\|\te^+\|_{L^2_k}-\frac14 \|\te\|_{L^2_k}-C'e^{-\delta_M t})
\\
\frac{\p}{\p t}\|\te^-\|_{L^2_k}&\le
-\delta(\|\te^-\|_{L^2_k}-\frac14 \|\te\|_{L^2_k}-C'e^{-\delta_M t})
\end{align}
\end{subequations}
for $t\in [T_{0,k},\infty)$. The arguments for the two inequalities are
very similar and we show (\ref{claim.ineq}b). We have

\begin{multline*}
-\frac12 \frac{\p}{\p t}\|\te^-\|_{L^2_k}^2=
-\langle\te^-,\frac{\p}{\p t}\te^-\rangle_k=
-\langle\te^-,\frac{\p}{\p t}\te\rangle_k=\\
\langle\te^-,*_Y d_Y\tilde\eta_3 - Q(\tilde\eta_3)
-\tilde\wp(t,\te)\rangle_k\ge
\delta \|\te^-\|_{L^2_k}^2-\frac14 \|\te\|_{L^2_k}\,\|\te^-\|_{L^2_k}
-C'e^{-\delta_M t}\,\|\te^-\|_{L^2}
\end{multline*}
where $\langle\cdot\, ,\cdot\rangle_k$ denotes the $L^2_k$ inner product
on~$Y$. This implies (\ref{claim.ineq}b) for all $t$ where $\te^+(t)\neq 0$
and, by continuity, on the closure of the set of such~$t$. The complement
is an open set where $\te^+(t)\equiv 0$ and the inequality is obvious.

We next check that $\|\te^+\|_{L^2_k}<\|\te^-\|_{L^2_k}$ for each sufficiently
large~$t$. Subtracting (\ref{claim.ineq}b) from (\ref{claim.ineq}a)
yields
$$
\frac{\p}{\p t}(\|\te^+\|_{L^2_k}-\|\te^-\|_{L^2_k})\ge
\frac\delta{2}\|\te\|_{L^2_k},
$$
so $\|\te^+\|_{L^2_k}-\|\te^-\|_{L^2_k}$ is non-decreasing in~$t$.
If our claim is false, then $\|\te^+\|_{L^2_k}-\|\te^-\|_{L^2_k}\ge 0$ on
some semi-infinite interval $t>T_0$. But then
$$
\frac{\p}{\p t}\|\te^+\|_{L^2_k}\ge
\delta(\|\te^+\|_{L^2_k}-\frac14 \|\te\|_{L^2_k})\ge
\delta(\|\te^+\|_{L^2_k}-\frac12 \|\te^+\|_{L^2_k})
$$
and $\|\te^+\|_{L^2_k}$ is exponentially increasing, a contradiction.
Finally,
$$
\frac{\p}{\p t}\|\te^-\|_{L^2_k}\le
-\delta(\|\te^-\|_{L^2_k}+\frac14\delta \|\te\|_{L^2_k}\le
-\delta\|\te^-\|_{L^2_k}+\frac12\delta \|\te^-\|_{L^2_k}
$$
implies the required $L^2_k$ estimate of $\te|_{\{t\}\times Y}$.
\end{pf}

\section{A class of asymptotically cylindrical Calabi--Yau
orbifolds}
\label{orbifolds}

In order to construct examples of asymptotically cylindrical
holonomy-$\Spin(7)$ manifolds by application of Theorem~\ref{exist.spin7},
we use complex projective 4-dimensional orbifolds equipped with
asymptotically cylindrical Calabi--Yau metrics.

Recall that the definition of an $n$-dimensional
orbifold, $V$ say, is analogous to that of a smooth manifold, but a
neighbourhood $U_p$ of a point $p$ is homeomorphic to $\RE^n/\Ga_p$,
where $a_p:\Ga_p \times \RE^n \to \RE^n$ is an effective action of a
finite group $\Gamma_p$ (the `local isotropy group') on~$\RE^n$. The
covering map $\RE^n \to U_p$ is called a {\em local uniformizing chart}
centred at $p$ and the coordinates on the domain $\RE^n$ are {\em local
uniformizing coordinates}. (Cf.~\cite{baily}.)

If $\Ga_p= \{1\}$ then $p$ is a {\em smooth point} of $V$, otherwise $p$ is
a {\em  singular} point. The set of all smooth points of $V$ is a manifold
and will be denoted $V^*$. We shall be interested in orbifolds whose
singular points are {\em isolated}, with each $\gamma\neq 1$ in $\Ga_p$
only fixing $0$ for each point~$p$.

The concepts of smooth differential forms and, more generally, sections of
vector bundles associated to $TV$ admit natural generalizations to
orbifolds: in a neighbourhood of singular point these are defined by the
respective $G$-invariant objects in uniformizing coordinates.

A Riemannian metric $g$ on $V^*$ is called
a {\em smooth orbifold metric}, and $(V,g)$ is called a Riemannian
orbifold, if the pull-back of $g$ to a local uniformizing chart
extends smoothly to $0\in \RE^n$ and $a_p$ acts by isometries. In
particular $a_p$ gives a representation of $\Ga_p$ in $O(n)$, the
orthogonal group of $T_0\RE^n$.
The holonomy group of a Riemannian orbifold $V$ is defined as the holonomy
group of its smooth part, the Riemannian manifold~$V^*$.

In a similar manner one defines $n$-dimensional complex orbifolds and
K\"ahler orbifold metrics, with each $\Gamma_p$ now assumed to be a
subgroup of $U(n)$. On the other hand, the canonical bundle of a complex
orbifold will be a well-defined orbifold complex line bundle if 
$\Ga_p\subset SL(n,\CX)$ for each~$p$ and then it makes sense to speak of
the anticanonical Weil divisors.

Our construction of asymptotically cylindrical $\Spin(7)$-manifolds will
start from asymptotically cylindrical Calabi--Yau complex 4-folds
obtained using the following result.
        \begin{thrm}\label{bCY}
Let $\oW$ be a compact K\"ahler 4-dimensional orbifold with K\"ahler
form~$\omega_{\oW}$ and with the singular locus of~$\oW$ consisting of
finitely many points with local isotropy groups contained in~$SU(4)$.
Suppose that $D\in|-K_{\oW}|$ is a smooth 3-fold in the anticanonical
class containing no singular points of~$\oW$ and such that the normal
bundle $N_{D/\oW}$ is holomorphically trivial. Suppose further that $D$ 
and the smooth locus $\oW^*\subset\oW$ and $\oW^*\setminus D$ are
simply-connected.

Then $W=\oW\setminus D$ admits a complete Ricci-flat K\"ahler metric $g$ with
holonomy~$SU(4)$ and a non-vanishing holomorphic volume $(4,0)$-form.
These are exponentially asymptotic to the product cylindrical Ricci-flat
K\"ahler structure on $D\times S^1\times \RE_{>0}$, via a diffeomorphism
$U\minus D\simeq D\times S^1\times \RE_{>0}$, where $U$ is a neighbourhood
of $D$ in $\oW$ and $z=\exp(-t-i\vartheta)$, $t>0$, $\vartheta\in S^1$ extends
over~$D$ to give a holomorphic coordinate on~$U\subset\oW$ vanishing 
to order one precisely on~$D$.
The K\"ahler form $\omega_g$ of~$g$ may be written near $D$ as
	\begin{equation}\label{asymp}
\omega_g|_{U\minus D} = dt\we d\vartheta + \omega_D + d\psi.
	\end{equation}
and a holomorphic volume form (up to a constant factor) as
	\begin{equation}\label{h.vol}
\theta_g|_{U\minus D} = (dt+id\vartheta)\we\theta_D + d\Psi.
	\end{equation}
Here $\omega_D$ is the Ricci-flat K\"ahler metric on~$D$ in the cohomology
class $[\omega_{\oW}|_D]$ and $\theta_D$ is a non-vanishing holomorphic
$(3,0)$-form on~$D$ satisfying $3i\,\theta_D\we\bar\theta_D = 4\,\omega_D^3$.
The forms $\psi,\Psi$ decay at the rate $O(e^{-\lambda t})$, $t\to\infty$,
with all derivatives, for some $\lambda>0$ depending on $\omega_D$.

The metric $g$ is uniquely determined by the above properties.
        \end{thrm}

Theorem~\ref{bCY} in the case of non-singular $\oW$ is an adapted from
\cite[\S\S 2--3]{g2paper}, \cite[pp.~142--143]{gokova} and \cite[\S~4]{HHN}.
The changes required in the presence of isolated orbifold singularities
away from~$D$ are straightforward as the foundational results used in the
proof carry over to orbifolds \cite{baily,satake}.
Note also that anticanonical divisor $D$ admits Ricci-flat K\"ahler
metrics as $c_1(D)=0$ by the adjunction formula.

The simply-connected condition on~$D$ is used for showing that the
holonomy of $W$ is $SU(4)$ (cf.~\cite[Theorem 2.7]{g2paper}). The hypothesis
can be weakened to $b^1(D)=0$ as the cross-section $D\times S^1$ then has
$b^1=1$; this requires an orbifold version of Theorem~\ref{hol8d}.

Recall that an antiholomorphic involution of a complex orbifold is defined
as a diffeomorphism $\rho$ of the underlying real orbifold such that
$\rho^2=\mathrm{id}$ and $\rho^*(J)=-J$, where $J$ is an almost complex
structure.
	\begin{prop}\label{spin7orbi}
Assume the hypotheses and notation of Theorem~\ref{bCY}.
Let $\rho$ be an anti-holomorphic involution of $\oW$, such that $\rho$
fixes precisely the singular points in~$\oW$ and maps $D$ onto itself.

Then there exists an asymptotically cylindrical Calabi--Yau structure
$(\omega,\theta)$ on the orbifold~$W$ satisfying the assertions of
Theorem~\ref{bCY} and such that
$$
\rho^*(\omega)=-\omega
\text{ and }
\rho^*(\theta)=\bar\theta.
$$
Furthermore, the 4-form
$$
\Phi=\frac12 \omega\we\omega+\Re\theta
$$
induces an asymptotically cylindrical torsion-free Spin(7)-structure on~$W$
which is invariant under $\rho$.
	\end{prop}
\begin{pf}
Let $g'$ be the K\"ahler metric of $\omega_{\oW}$. We may assume without loss
that $g$ is invariant under $\rho$, by passing if necessary to
$\frac12 (g+\rho^*(g))$. Then the antiholomorphic involution of~$D$
defined by the restriction of~$\rho$ preserves the K\"ahler metric on~$D$
induced by~$g$. The argument in~\cite[Prop.~15.2.2]{joyce} applies with
only cosmetic changes to $D$ show that $\rho^*(\omega_D)=-\omega_D$ and
$\rho^*(\theta)=\bar\theta$.

Let $g$ be the asymptotically cylindrical Calabi--Yau metric on~$W$ given by
Theorem~\ref{bCY} with K\"ahler form asymptotic to $\omega_D$.
It is not difficult to check using \cite[Prop.~1.1]{gokova} that the
K\"ahler metric $\rho^*(g)$ is also asymptotic to $\omega_D$ in the 
sense of Theorem~\ref{bCY}. Then $\rho^*(g)=g$ by the uniqueness of~$g$.
\end{pf}
\begin{rmk}
Note also that $\rho$ induces an antiholomorphic involution on the
canonical bundles of $\oW$ and $D$, hence by adjunction on the conormal
bundle of~$D$. It can be checked that a local holomorphic coordinate $z$
vanishing to order 1 on~$D$ may be chosen so that $z\circ\rho=\bar z$.
Thus the action of $\rho$ on the cylindrical end of~$\oW$ is asymptotic to
the antiholomorphic involution $\rho_\infty$ on $(\RE\times S^1)\times D$
acting as $t+i\vartheta\mapsto t-i\vartheta$ on the first factor and as
$\rho|_D$ on the second factor. The asymptotic torsion-free
$\Spin(7)$-structure is
\begin{equation}\label{cyl}
\Phi_\infty=dt\we d\vartheta\we\omega_D+dt\we\Re\theta_D-
d\vartheta\we\Im\theta_D+\frac12\omega_D\we\omega_D,
\end{equation}
induced by the Calabi--Yau structure $(\omega_D,\theta_D)$ on the
3-fold~$D$ (cf.~\cite[Prop.~13.1.2]{joyce}), and is invariant
under~$\rho_\infty$. 
\end{rmk}

The next proposition will be useful for providing examples of orbifolds
satisfying the hypotheses of Theorem~\eqref{bCY}.
	\begin{prop}[{cf.~\cite[Prop.~6.42]{g2paper}}] \label{blowup}
Let $V$ be a compact K\"ahler 4-dimensional orbifold and suppose that
a 3-fold $D\subset V^*$ is a smooth anticanonical divisor $D\in|-K_V|$.
Let $\Sigma$ be a connected smooth surface in $D$ representing the
self-intersection $D\cdot D$ in the Chow ring of~$V$ and denote by
$\sigma:\tV\to V$ the blow-up of $V$ along~$\Sigma$.

Then the closure $\tD$ of $\sigma^{-1}(D\minus\Sigma)$ is a smooth
anticanonical divisor on $\tV$ and \mbox{$\tD\cdot\tD=0$}, so the normal
bundle of $\tilde{D}$ is holomorphically trivial. Furthermore,
$\sigma$ restricts to give an isomorphism $\tD\to D$ of complex
3-folds and this isomorphism identifies the K\"ahler metric
on $D\subset V$ with the restriction to $\tD$ of some K\"ahler metric
on~$\tV$.
	\end{prop}
        \begin{rmk}
If the orbifold $V$ satisfies the hypotheses of Proposition~\ref{blowup}
and both $V^*$ and $D$ are simply-connected, then $\tV\minus\tD$ is also
simply-connected because the exceptional divisor on~$\tV$, which is the
projective normal bundle of $\Sigma$, has a fibre $\CP^1$ meeting the
proper transform~$\tD$ in exactly one point.
	\end{rmk}
	\begin{rmk}\label{blow.div}
Proposition~\ref{blowup} is easily generalized to the case when $\Sigma$ is
a finite sequence of connected smooth surfaces $\Sigma_1,\ldots,\Sigma_m$.
These need not necessarily be distinct, i.e.\ the argument allows an
effective divisor on~$D$ with smooth components of multiplicity $>1$.
One then recursively blows up each $\Sigma_i$ and lifts the remaining
$\Sigma_{i+1},\ldots,\Sigma_m$ to the respective proper transform of~$D$.
Cf.~\cite[p.~199]{KL}.
	\end{rmk}

Assuming further that the singular locus of $V$ is finite with each local
isotropy group in $SU(4)$, we can apply Theorem~\ref{bCY} to
$\overline{W}=\tV$ and obtain on $W=\tV\minus\tD$ an asymptotically
cylindrical Ricci-flat K\"ahler metric~$g$ with holonomy $SU(4)$. Then $W$
is also an instance of asymptotically cylindrical $\Spin(7)$-manifold and
$g$ is induced by the $\Spin(7)$-structure defined via~\eqref{su4spin7}.

Now suppose that $V$ also has an antiholomorphic involution $\rho$ fixing
precisely the singular points of~$V$ and mapping each of $D$, $\Sigma$ onto
itself. Then $\rho$ lifts to an antiholomorphic involution $\trho$ of $\tV$
so that
\begin{equation}\label{trho}
\begin{CD}
\tV @>\trho>> \tV\\
@V\sigma VV	@VV\sigma V\\
V @>\rho>> V
\end{CD}
\end{equation}
is a commutative diagram. Further, $\trho$ defines by restriction an
antiholomorphic involution of~$\tD$. Therefore, by
Proposition~\ref{spin7orbi}, the asymptotically cylindrical
$\Spin(7)$-structure on $W=\tV\minus\tD$ can be taken to be
$\trho$-invariant and descends to the quotient
\begin{equation}\label{M.orb}
M_0=(\tV\minus\tD)/\trho,
\end{equation}
and $M_0$ is an asymptotically cylindrical $\Spin(7)$-orbifold. As
$\tV^*\minus D$ is simply-connected and $\trho$ acts freely the fundamental
group of the smooth locus $M_0^*\subset M_0$ is $\ZE_2$.

The $G_2$-structure $\phi$ `at infinity' on the cross-section
$Y_W=S^1\times D$ of~$W$ is induced by the Calabi--Yau structure on the
divisor,
$$
\phi=d\vartheta\we\omega_D + \Re\theta_D.
$$
The antiholomorphic involution $\rho$ acts freely on the 3-fold $D$
and as reflection on~$S^1$, thus freely of~$Y_W$ preserving the `product'
$G_2$-structure, $\rho^*\phi=\phi$. As $b^1(D)=0$ we obtain that the
cross-section $Y\cong (D\times S^1)/\rho$ of $M_0$ has $b^1(Y)=0$ and
infinite $\pi_1(Y)$. The holonomy of the metric $g(\phi)$ induced by the
$G_2$-structure on~$Y$ is a semi-direct product of $\ZE_2$ and $SU(3)$.

We next impose a condition on the singularities of $V$. As $D$ contains no
singularities of~$V$ the fixed points of~$\trho$ are preimages of the
fixed points of~$\rho$ and the respective singularities of $V$ and $\tV$
are isomorphic. The type of singularities and the method of resolution of
these are adapted from \cite[Chap.~15]{joyce}, the only difference being
that we shall preform the resolution on asymptotically cylindrical, rather
than compact Calabi--Yau 4-folds.

We require that the orbifold singularities of $V$ have local isotropy group
$\ZE_4$, so that the action of a generator of $\ZE_4$ is given in
uniformizing coordinates by
\begin{equation}\label{cyclic}
(z_1,z_2,z_3,z_4)\mapsto (iz_1,iz_2,iz_3,iz_4).
\end{equation}
With this arrangement, $\ZE_4$ acts as a subgroup of $SL(4,\CX)$, so the
canonical bundle of $V$ (and of $\tV$) is a well-defined `orbifold line
bundle' and $c_1(V)$ is a well-defined class in $H^2(V,\ZE)$.

The resulting singularities of $M_0$ have non-abelian local isotropy group
$G$ on two generators coming from $\ZE_4$ and $\rho$. In fact, the
singularities of~$M_0$ are of the same type; the action of the generators of
each is isomorphic to $G\subset\Spin(7)$ generated by the unitary maps
\eqref{cyclic} and
\begin{equation}\label{involution}
(z_1,z_2,z_3,z_4)\mapsto (\bar z_2,-\bar z_1,\bar z_4,-\bar z_3),
\end{equation}
with the real and imaginary parts of $\p/\p z_1,\ldots,\p/\p z_4$
corresponding to an orthonormal frame at $p_i$ in uniformizing real
coordinates (cf. \cite[Prop.~15.2.3]{joyce}).

The orbifold singularities modelled on $\RE^8/G$ can be resolved using a
method given in~\cite{joyce}, by forming at each singular
point $p_i$ ($i=1,\ldots,k$) a generalized connected sum of $M_0$ and an
asymptotically locally Euclidean (ALE) $\Spin(7)$-manifold asymptotic to
$\RE^8/G$. There are two suitable choices of ALE $\Spin(7)$-manifolds
$(X_1,\Phi_1)$, $(X_2,\Phi_2)$ with $\pi_1(X_i)=\ZE_2$ for $i=1,2$
(\cite[\S 15.1]{joyce}). The respective $\Spin(7)$-compatible generalized
connected sums of $M_0\# X_i$ correspond to different identifications on the
cross-section of the neck $S^7/G$ and the resulting spaces are topologically
distinct. 

\begin{prop}[{cf. \cite[Prop. 15.2.9 and 15.2.13]{joyce}}]\label{esti}
The asymptotically cylindrical $\Spin(7)$-orbifold defined in~\ref{M.orb}
admits a simply-connected resolution of singularities\linebreak
$M_0\# (X^{(1)}\sqcup\ldots\sqcup X^{(k)})$ with $X^{(j)}$ are ALE
$\Spin(7)$-manifolds asymptotic to $\RE^8/G$.

The smooth 8-manifold $M$ admits a 1-parameter family of asymptotically
cylindrical\linebreak $\Spin(7)$-structures $\Phi^s$, $0<s<\eps$, which are
torsion-free on the cylindrical end. There are positive constants
$\lambda,\mu,\nu$, such that
\begin{enumerate}
\itemsep0em
\item
the injectivity radius $\delta(g)$ satisfies $\delta(g)\ge\mu s$,
\item
the Riemann curvature $R(g)$ satisfies $\|R(g)\|_{C^0}\le\nu s^{-2}$,
\item
$\|\psi\|_{L^2}\le\lambda s^{24/5}$
and
$\|d\psi\|_{L^{10}}\le\lambda s^{36/25}$,
\end{enumerate}
for all $0<s<\eps$, with the norms defined using the metric $g(\Phi^s)$.
\end{prop}

Parts (a),(b),(c) of Proposition~\ref{esti} imply the hypotheses (a),(b)
and \eqref{psi} of the existence Theorem~\ref{exist.spin7} and taking a
sufficiently small $s>0$ we obtain an asymptotically cylindrical
torsion-free structure $\Phi$ on~$M$.
Since the cross-section of~$M$ has $b^1=0$ and $M$ has single end, the
induced metric  $g(\Phi)$ has holonomy $\Spin(7)$ by Theorem~\ref{hol8d}(a).

Putting together the conditions imposed on $V$ we make the following.

\begin{defi}\label{config}
We say that $(V,D,\Sigma,\rho)$ is an {\em admissible orbifold configuration}
if the following properties hold.

(i) $V$ is a compact K\"ahler 4-dimensional orbifold. The singular locus
of~$V$ is non-empty and consists of finitely many isolated singular
points $\{p_1,\ldots,p_k\}$. Each singularity is isomorphic
to~\eqref{cyclic} with the local isotropy group~$\ZE_4$.

(ii) The smooth locus $V^*\subset V$ is simply-connected and there exists a
simply-connected smooth 3-fold $D\subset V^*$ in the anticanonical class
$D\in|-K_V|$ and $V^*\minus D$ is simply-connected.
The self-intersection $D\cdot D$, in the sense of the Chow ring,
is represented by an effective divisor $\Sigma=\sum_i n_i\Sigma_i$, where
each $\Sigma_i$ is a smooth complex surface in~$D$.

(iii)
There is an antiholomorphic involution $\rho$ of $V$ fixing 
the singular locus and no other points in~$V$ and such that
$\rho(D)=D$ and $\rho(\Sigma_i)=\Sigma_i$, for each $i$.

We call a 4-orbifold $V$ {\em admissible} if $V$ is part of some admissible
orbifold configuration.
\end{defi}

Summarizing the work in this section, we obtain.

\begin{thrm}\label{construct}
Suppose that $(V,D,\Sigma,\rho)$ is an admissible orbifold configuration
and let $\tV$ be obtained from $V$ by successive blow-up of the components
of $\Sigma$ according to multiplicities. Let $\tD$ be the proper transform
of~$D$ and $\trho$ the lift of $\rho$ to~$\tV$.

Then the real 8-dimensional orbifold $M_0=(\tV\minus\tD)/\trho$ admits a
simply-connected resolution of singularities $M$, so that there is an
asymptotically cylindrical metric on $M$ with holonomy equal to $\Spin(7)$.
The cross-section of $M$ `at infinity' is isometric to $(D\times S^1)/\rho$
and has holonomy $\ZE_2 \ltimes SU(3)$, where $D$ is taken with the
Ricci-flat K\"ahler metric determined by the restriction of the K\"ahler
class of~$V$.
\end{thrm}

We shall refer to $M$ in this theorem as the $\Spin(7)$-manifold
{\em constructed from $(V,D,\Sigma,\rho)$}.

\section{Examples of the construction} 
\label{examples}

By the results of the previous section, in order to construct
asymptotically cylindrical 8-manifolds with holonomy $\Spin(7)$ it suffices
to find admissible orbifolds. In this section, we give simple examples
using the geometry of weighted complex projective spaces.

For a set of positive integers $a=\{a_0,\ldots,a_n\}$ with highest common
factor~1, the complex {\em weighted projective space} $\CP^n_a$ is
defined as the quotient of $\CX^{n+1}\minus\{0\}$ by the equivalence relation
$$
(z_0,\ldots,z_n)\sim (u^{a_0}z_0,\ldots,u^{a_n}z_n),
\text{ for each }u\in\CX^*
$$
(when all $a_i=1$ this gives the usual complex projective space).
We shall denote by $[z_0,\ldots,z_n]$ the equivalence classes. Each $\CPA$
is a projective complex orbifold and admits K\"ahler metrics. The smooth
locus $(\CPA)^*$ is simply-connected~\cite{dimca85}
and the rational cohomology ring of $\CPA$ is $\QU[x]/x^{n+1}$ generated,
like in the case of the usual complex projective space, by an element $x$
of degree~2 \cite[Prop.~B13]{dimca-book}.

A number of standard properties of varieties in $\CP^n$ extends to weighted
projective spaces. For the proof and further details of results used below,
we refer to \cite{idolga}, \cite[Appendix~B]{dimca-book} and~\cite{IF}.

If a polynomial $f(z_0,\ldots,z_n)$ satisfies a {\em weighted homogeneous}
condition of degree~$d$
$$
f(u^{a_0}z_0,\ldots,u^{a_n}z_n)=u^df(z_0,\ldots,z_n)
\text{ for each }u,z_0,\ldots,z_n\in\CX,
$$
then the zero locus of $f$ is a well-defined hypersurface $S$ in
$\CPA$ and $d$ is the degree of~$S$.
Hypersurfaces and complete intersections of (complex) dimension $\ge 2$ in
a weighted projective space are simply-connected.

A hypersurface $S$ is called a {\em linear cone} if the defining weighted
homogeneous polynomial may be written as  $f=z_i+g$ for some $i$. In this
case, $\{f=0\}$ is isomorphic to an $(n-1)$-dimensional weighted projective
space with weights $a_0,\ldots,\hat a_i,\ldots,a_n$. Here and later in this
section the \^\ notation means this element is omitted.

If there are no points $(z_0,\ldots,z_n)\neq 0$ such that
$f(z_0,\ldots,z_n)=0$ and $\p f/\p z_j(z_0,\ldots,z_n)=0$, then 
$S$ has no singularities outside the singular locus of $\CP^n_a$ and
any singularity of $S$ is of orbifold type with a cyclic group
\cite[Theorem~3.1.6]{idolga}.

Each $\CPA$ may be considered as an algebraic variety and the definition of
sheaves $\OH(m)$ extends to weighted projective spaces. Furthermore,
the anticanonical class of a hypersurface or complete intersection, $Z$
say, is determined by application of the adjunction formula if one assumes 
that the preimage of $Z$ is a smooth cone in $\CX^{n+1}\minus\{0\}$ (the
{\em quasismooth} condition) and $Z$ is {\em well-formed in $\CPA$}. For
hypersurfaces of degree~$d$, the well-formed condition means

$\bullet$
$\gcd (\ldots,\hat{a}_i,\ldots)=1$ for all $i$;

$\bullet$
$\hcf (\ldots,\hat{a}_i,\ldots,\hat{a}_j,\ldots)$ divides $d$, for
all $i\neq j$.
\\
For a complete intersection of two hypersurfaces, of degrees $d_1$ and
$d_2$, the condition becomes 

$\bullet$
$\hcf (\ldots,\hat{a}_i,\ldots)=1$ for all $i$;

$\bullet$
$\hcf (\ldots,\hat{a}_i,\ldots,\hat{a}_j,\ldots)$ divides both $d_1$ and
$d_2$, for all $i\neq j$;

$\bullet$
$\hcf (\ldots,\hat{a}_i,\ldots,\hat{a}_j,\ldots,\hat{a}_k,\ldots)$ divides
one or both of $d_1$ and $d_2$, for all distinct $i,j,k$.
\\
A general definition requires more notation and can be found
in~\cite[\S 6]{IF}.

It can be shown that the sum of linear cones $\{z_i=0\}$, for all
$i=0,\ldots,n$, defines an anticanonical (Weil) divisor on $\CPA$.
The anticanonical sheaf on a well-formed complete intersection $S$ of
$k$ hypersurfaces of degree $d_1,\ldots,d_k$ in $\CP^n_a$ then is given by
\begin{equation}\label{c1}
-\mathcal K_S  = \OH(\sum_{j=0}^n a_j - \sum_{i=0}^k d_i).
\end{equation}
Thus $\sum_{i=0}^k d_i < \sum_{j=0}^n a_j$ is a necessary condition for $S$ to
have effective anticanonical divisors.

There is a version of Lefschetz hyperplane theorem.
\begin{thrm}[{\cite[Theorem~B22]{dimca-book}}]
If $Z$ is a hypersurface or complete intersection in~$\CPA$, then the
homomorphism $H^r(\CPA;\QU)\to H^r(Z,\QU)$ induced by the inclusion is an
isomorphism for $r<\dim_{\CX}Z$ and is injective for $r=\dim_{\CX}Z$.
\end{thrm}

\subsection{Considerations of cohomology}
\label{de.Rham}

We next determine the de Rham cohomology of $\Spin(7)$-manifolds obtainable
via the procedure of \S\ref{orbifolds} when $V$ is an orbifold in a
weighted projective space. Recall from \S\ref{as.cyl} and Theorem~\ref{moduli}
the definitions of $b^r_0$ and $b^4_\pm$.

\begin{prop}\label{betti}
Assume $(V,D,\Sigma,\rho)$ is an admissible orbifold configuration and
a hypersurface or complete intersection $V\subset\CP^n_a$ is quasismooth
and well-formed. Then the Betti numbers of the $\Spin(7)$-manifold $M$ with
cross-section~$Y$ constructed from $(V,D,\Sigma,\rho)$ are
\begin{align*}
b^1(Y)=b^2(Y)=0, &\qquad b^3(Y)=2+h^{2,1}(D),\\
b^1(M)=b^1_c(M)=b^2(M&)=b^2_c(M)=b^3(M)=b^3_c(M)=0,\\
b^4(M)=b^4_0(M)+b^3(Y), &\qquad
b^4_0(M)=\frac12(\chi(\Sigma)+\chi(V)+3k)-4,
\end{align*}
where $k$ is the number of singular points in~$V$.
\end{prop}
\begin{rmk}
Every asymptotically cylindrical $\Spin(7)$-manifold $M$ has $b^4(M)>0$ as
the closed admissible 4-form $\Phi$ is non-decaying bounded harmonic and
defines a non-trivial cohomology class on a cross-section of cylindrical
end. The cross-section $Y$ is a compact $G_2$-manifold and $b^3(Y)>0$ as
there is a non-trivial harmonic 3-form $\phi$ inducing the $G_2$-structure.
Thus no further Betti numbers of $M$ and $Y$ in Proposition~\ref{betti} can
vanish.
\end{rmk}
\begin{pf}[Proof of Proposition~\ref{betti}]
The vanishing of first Betti numbers follows at once
as $M$ is simply-connected and admits holonomy $\Spin(7)$ metrics.
We have $H^j(\CP^n_a)\cong H^j(\CP^n)$ for each $j$ and then 
$b^1(V)=b^1(D)=b^1(\Sigma)=b^1(\CP^n_a)=0$ and
$b^2(V)=b^2(D)=b^2(\CP^n_a)=1$ and $b^3(V)=b^3(\CP^n_a)=0$
by application of the Lefschetz hyperplane theorem.
The cohomology $H^*(V/\rho)$ is the $\rho$-invariant part 
$H^*(V)^\rho\subset H^*(V)$ and the cohomology in degree 2 is generated by
the K\"ahler form on which 
the antiholomorphic involution acts as $-1$. Thus
$b^2(V)^\rho=b^2(D)^\rho=0$.

The cohomology of the blow-up may be written as
$H^r(\tV)\cong H^r(V)\oplus (H^r(E)/\sigma^* H^r(\Sigma))$, where $E$ is
the exceptional divisor \cite[p.~605--608]{GH}. The cohomology of $E$ is
generated as a cup-product algebra over $H^*(\Sigma)$ by the Chern class
$\zeta=c_1([E])|_E$ with a relation
$\zeta^2-c_1(N_{\Sigma/M}) \zeta+c_2(N_{\Sigma/M})=0$, where $N_{\Sigma/M})$
is the normal bundle. The pull-back $\sigma^*$ is injective on
$H^*(\Sigma)$, so $H^j(E)/\sigma^*(\Sigma)$ is generated by 
$\zeta\we \sigma^*H^{j-2}(\Sigma)$ and the blow-up adds $b^{j-2}(\Sigma)$
to the $j$-th Betti number. Now ${\trho\,}^*\zeta=-\zeta$ as $c_1([E])$ is
Poincar\'e dual to the cycle $[E]$ defined by complex 3-dimensional
submanifold on which $\trho_*$ changes the orientation. We deduce that
$b^j(\tV)^{\trho}=b^j(V)^{\trho}+b^{j-2}(\Sigma)^{-\trho}$ and
\begin{gather*}
b^1(\tV)^{\trho}=b^2(\tV)^{\trho}=b^3(\tV)^{\trho}=0,
\qquad
b^2(\tV)^{-\trho}=h^{1,1}(\tV)^{-\trho}=2,
\\
b^4(\tV)^{\trho}=b^4(V)^\rho + b^2(\Sigma)^{-\rho}
\end{gather*}
The Betti numbers in middle dimension are now determined by the Euler
characteristic, $b^2(\Sigma)^{-\rho}=b^2(\Sigma)-b^2(\Sigma)^\rho=
\frac12\chi(\Sigma)$ as $\rho$ acts freely on~$\Sigma$ and $b^1(\Sigma)=0$.
Similarly, $b^4(V)^\rho=\frac12(\chi(V)+k)-2$.

For the Calabi--Yau 3-fold $D$, we have $b^3(D)^\rho=b^3(D)^{-\rho}=\frac12
b^3(D)=1+h^{2,1}(D)$ as $\rho$ interchanges $H^{i,j}(D)$ and $H^{j,i}(D)$.
Then for $Y\cong (D\times S^1)/\rho$, recalling that $b^1(Y)=0$ and that
$\rho$ acts on $S^1$ as reflection, we obtain
$$
b^2(Y)=b^2(D)^\rho=0,\qquad  b^3(Y)=2+h^{2,1}(D)
$$
as $\rho^*\omega_D=-\omega_D$ from Proposition~\ref{spin7orbi}
and~\eqref{cyl}.

The resolution of each singularity of $\tV/\trho$ is topologically a
generalized connected sum $\oM$ of $\tV/\trho$ and an ALE $\Spin(7)$ manifold
$X_i$ and cross-section of the neck is a spherical space-form $S^7/G$.
The manifolds $X_i$ are taken from \cite[\S 15.1]{joyce} and have
$b^1(X_i)=b^2(X_i)=b^3(X_i)=0$, $b^4(X_i)=b^4_-(X_i)=1$, so the resolved
compact smooth 8-manifold $\oM$ has $b^j(\oM)=b^j(\tV)^{\trho}=0$,
for $j=1,2,3$, and
$$
b^4(\oM)=b^4(\tV)^{\trho}+k= b^4(V)^\rho+b^2(\Sigma)^{-\rho}+k.
$$

To determine the cohomology of asymptotically cylindrical $M$ with
cross-section $Y$ we use the Mayer--Vietoris theorem applied to
$\oM=M\cup U$. Here $U$ is a tubular neighbourhood of $\tD/\trho$ and
$M\cap U$ retracts to~$Y$. Firstly, $b^1(M)=0$ as $M$ is 
simply-connected and also $b^2(M)=b^2(\oM)=0$ as $b^1(Y)=b^2(Y)=0$. In the
part of the Mayer--Vietoris exact sequence
\begin{equation}\label{deg3}
0\to H^3(\oM)\to H^3(M)\oplus H^3(D)^\rho\to H^3(D)^\rho\oplus\RE[\phi]
\end{equation}
the last homomorphism maps onto $H^3(D)^\rho$ as the image contains
$H^3(D)^\rho$ and is orthogonal to the $G_2$ 3-form $\phi$ on~$Y$ by
Proposition~\ref{3-forms}. We find that $b^3(M)=b^3(\oM)=0$. Then
$b^2_c(M)=0$, $b^3_c(M)=0$ and $b^4_0(M)=b^4_c(M)-b^3(Y)=b^4(M)-b^3(Y)$
from the long exact sequence \eqref{exact} for de Rham cohomology with
compact support and the Poincar\'e duality. Extracting a further part of
the Mayer--Vietoris exact sequence
\begin{equation}\label{deg4}
0\to \RE[\phi]\to H^4(\oM)\to H^4(M)\oplus\RE[\omega_D^2]\to
\RE [\omega_D^2]\oplus H^3(D)^{-\rho}\to 0, 
\end{equation}
we deduce $b^4(M)=b^4(\oM)+b^3(D)^{-\rho}-1=b^4(\oM)+b^3(Y)-2$. Then
$b^4_0(M)=b^4(\oM)-2=b^4(V)^\rho+b^2(\Sigma)^{-\rho}+k-2$ and the required
formula for $b^4_0(M)$ follows by substitution of the expressions for
$b^4(V)^\rho$ and $b^2(\Sigma)^{-\rho}$ found above.
Proposition~\ref{betti} is proved.
\end{pf}

The only non-trivial multiplicative structure on the de Rham cohomology
of~$M$ is the intersection form on $H^4_c(M)$. A maximal subspace on which
the latter form is non-degenerate is isomorphic to $H^4_0(M)$.

\begin{prop}\label{sig}
In the situation of Proposition~\ref{betti},
the signature of the intersection form on $H^4_0(M)$ is
\begin{align*}
&b^4_+(M)=h^{2,2}(V)^{\rho}+h^{1,1}(\Sigma)^{-\trho}-2,
\\
&b^4_-(M)=h^{3,1}(V)+h^{2,0}(\Sigma)+k.
\end{align*}
\end{prop}
Note that $b^4_+(M)+b^4_-(M)=b^4_0(M)$ recovers the expression for
$b^4_0(M)$ in Proposition~\ref{betti} because
$h^{2,2}(V)^{\rho}+h^{3,1}(V)=b^4(V)^\rho$ and 
$h^{1,1}(\Sigma)^{-\trho}+h^{2,0}(\Sigma)=\frac12\chi(\Sigma)$.
\begin{pf}
The signature of K\"ahler orbifold $\tV$ is
$b^4_+(\tV)^{\trho}=h^{2,2}(\tV)^{\trho}-h^{1,1}(\tV)^{-\trho}+1$, \
$b^4_-(\tV)^{\trho}=h^{3,1}(\tV)+h^{1,1}(\tV)^{-\trho}-1$, determined from
the Hodge--Riemann bilinear relations.
The ALE $\Spin(7)$-manifolds used in the resolution of singularities of
$\tV$ have $b^4_-(X_i)=b^4(X_i)=1$ and it follows from the argument of
Proposition~\ref{betti} that
$b^4_+(\oM)=h^{2,2}(V)^{\rho}+h^{1,1}(\Sigma)^{-\trho}-1$ and
$b^4_-(\oM)=h^{3,1}(V)+h^{2,0}(\Sigma)+k+1$.

Comparing the exact sequence \eqref{deg4} and
$$
0\to \RE [\phi]\oplus H^3(D)^{\rho}\to H^4_c(M)\to H^4(M)\to 
\RE [\omega_D^2]\oplus H^3(D)^{-\rho}\to 0
$$
we find that the image of $H^4(\oM)$ in $H^4(M)$ contains $H^4_0(M)$ as a
codimension 1 hyperplane. Its complement is spanned by a class $\xi'$ such
that $\xi'|_{\{t_0\}\times Y}=[\omega_D^2]$ and $\xi'$ may be taken orthogonal
to $H^4_0(M)$ with respect to the cup product. The kernel of
$H^4(\oM)\to H^4(M)$ is spanned by $\xi''=[d\eta(t)\phi]$ 
where $\eta(t)$ is a cut-off function with derivative supported on the
cylindrical end of~$M$.

On the other hand, $H^4_0(M)$ is the subspace of classes in $H^4(M)$ that
are represented by closed 4-forms with compact support in~$M\subset\oM$.
We deduce that the intersection form on $H^4_0(M)$ may be identified with the
intersection form on the orthogonal complement in $H^4(\oM)$ of the 2-plane
spanned by $\xi',\xi''$. As $\xi'\cupprod\xi''\neq 0$ and $\xi'\cupprod\xi'=0$
this latter 2-plane has signature $(1,1)$ and $H^4_0(M)$ has the complementary
signature in $H^4(\oM)$.
\end{pf}

\subsection{An admissible weighted projective space}
\label{wps}

Let $V=\CP^4_{1,1,1,1,4}$. This orbifold has unique
singular point \mbox{$p_0=[0,0,0,0,1]$} and the local isotropy group is
$\ZE_4$ with the required action~\eqref{cyclic}. The anticanonical sheaf of
$\CP^4_{1,1,1,1,4}$ is $\OH(8)$ and an anticanonical Weil divisor $D$
may be given by the vanishing of weighted degree~8 homogeneous polynomial
$$
f(z)=z_0^8+z_1^8+z_2^8+z_3^8+z_4^2.
$$
As $df(z)$ does not vanish for $z\neq 0$ and $D$ does not contain $p_0$
the hypersurface $D$ is smooth. In fact, $\CP^4_{1,1,1,1,4}$ is
the unique non-smooth weighted projective space with only $\ZE_4$-singularities
and having a smooth Calabi--Yau 3-fold in the anticanonical class,
cf.~\cite[Theorem~14.3]{IF}.

The hypersurface $D$ is well-formed and the Euler characteristic of~$D$ is
determined via the Chern class (\cite[p.~393]{candelas}) from the
adjunction formula, $c(D)= (1+x)^4(1+4x)(1+8x)^{-1}$, whence $\chi(D)=
\langle -148x^3,[D]\rangle=(-148/4)\cdot 8=-296$ and $h^{2,1}(D)=149$
(noting also that $\CP^4\to\CP^4_{1,1,1,1,4}$ is a 4-to-1 branched cover).

Choosing another anticanonical divisor $D'$ defined by the vanishing of
$$
f_1(z)=z_0^8-z_1^8+2z_2^8-2z_3^8+iz_4^2
$$
we obtain $\Sigma_{8,8}=D\cap D'$ a well-formed smooth complete intersection
surface. From $c(\Sigma_{8,8})=(1+x)^4(1+4x)(1+8x)^{-2}$, the Euler
characteristic is
$\chi(\Sigma_{8,8})=\langle 86x^2,[\Sigma_{8,8}]\rangle=86\cdot(1/4)\cdot 64=1376$.
From Noether's formula for the holomorphic Euler characteristic 
$\chi(\OH_{\Sigma_{8,8}})$ we determine the geometric genus
$h^{2,0}(\Sigma_{8,8})=199$.

The anti-holomorphic involution of $\CP^4_{1,1,1,1,4}$ given by
$$
\rho_1: [z_0,\ldots,z_4]\mapsto
[\bar{z}_1,-\bar{z}_0,\bar{z}_3,-\bar{z}_2,\bar{z}_4]
$$
fixes only $p_0$ and preserves each of $D$, $D'$ and $\Sigma_{8,8}$. Thus all the
requirements of Definition~\ref{config} are satisfied and
$(\CP^4_{1,1,1,1,4},D,\Sigma_{8,8},\rho_1)$ is an admissible orbifold configuration.

The 8-manifold $M_1$ with asymptotically cylindrical, holonomy $\Spin(7)$
metrics arising from this configuration by Theorem~\ref{construct} 
has $b^4_0(M_1)= \frac12\chi(\Sigma_{8,8})=688$, $b^4(M_1)=839$
and the cross-section has $b^3(Y)=151$ by Proposition~\ref{betti}.
From Theorem~\ref{moduli} and Proposition~\ref{sig}, we obtain
$b^4_-(M_1)=1+p_g(\Sigma_{8,8})=200$ and there is a
352-dimensional moduli space for asymptotically cylindrical
$\Spin(7)$-structures on~$M_1$.

\subsection{An admissible hypersurface}

Consider a hypersurface of weighted degree~8 in a 5-dimensional weighted
projective space 
\begin{equation}\label{hsurf}
V=\{f_8(z)=z_0^{8}+z_1^{8}+z_2^{8}+z_3^{8}+z_4^2+z_5^2=0\} \subset\CPW.
\end{equation}
The singular locus of $\CPW$ is a copy of the projective line $\CP^1$
given by $\{z_0=z_1=z_2=z_3=0\}$. Every singularity of $\CPW$ is locally
modelled on $(\CX^4/\ZE^4)\times\CX$. 
The singular locus of $V$ consists of two isolated points
$p_\pm=[0,0,0,0,i,\pm 1]$ which are orbifold singularities of
type~\eqref{cyclic}. The smooth locus $V^*$ is simply-connected by the
results in~\cite[Prop.~6 and Cor.~9]{dimca85}.

The hypersurface $V$ is well-formed and $-K_V=\OH(4)|_V$ by the
formula~\eqref{c1} and an anticanonical divisor $D$ on~$V$ may be given by
the intersection of $V$ with a linear cone
$C=\{z_4+z_5=0\}\cong\CP^4_{1,1,1,1,4}$. The unique singular point of~$C$
is $p_0=[0,0,0,0,1,-1]$ and $D=C\cap V$ contains no singular points of $C$
or $V$. It is not difficult to see that the tangent spaces of $C$ or $V$
are transverse at each point of~$D$ and the 3-fold $D$ is a smooth complete
intersection, isomorphic to the Calabi--Yau 3-fold in \S\ref{wps}.

The self-intersection of $D$ in the present case is realized by
$\Sigma_8=V\cap C\cap C'$, with $C'=\{z_4-z_5=0\}$ and is isomorphic to the
degree 8 Fermat surface $\Sigma_8=\{z_0^8+z_1^8+z_2^8+z_3^8=0\}$ in $\CP^3$. 
Its Euler characteristic is $\chi(\Sigma_8)=304$, determined again via the
Chern class. Further, using Noether's formula for the holomorphic Euler
characteristic, we find $h^{2,0}(\Sigma_8)=35$.

Each of $V,C,C'$, hence also $D$ and $\Sigma_8$, are invariant under an
antiholomorphic involution of $\CPW$ defined by 
$$
\rho_2:[z_0,\ldots,z_5]\mapsto
[\bar{z}_1,-\bar{z}_0,\bar{z}_3,-\bar{z}_2,\bar{z}_5,\bar{z}_4]
$$
and the fixed point set of $\rho_2$ on $V$ is precisely the singular locus
$p_\pm$. The configuration $(V,D,\Sigma_8,\rho_2)$ is admissible.

We calculate $\chi(V)=306$ via a recursive method given in
\cite[Ex.~15.3.4, cf.~also \S 15.4]{joyce}, using the nested sequence
of orbifolds $V_j=V\cap\{z_{j+1}=\ldots=z_5=0\}$ and the
branched covers $[z_0,\ldots,z_j]\in V_j\to [z_0,\ldots,z_{j-1}]
\in\CPW\cap\{z_{j}=\ldots=z_5=0\}$, for $j=0,1,\ldots,5$.
(Note that as $V$ has orbifold singularities the method using the Chern
class computes a different quantity known as `orbifold Euler
characteristic' which in general is {\em not} $\sum_j (-1)^jb^j(V)$ and in
the present case is a not an integer.)
The Hodge numbers of $V$ can be computed by Steenbrink's
result~\cite[Theorem~7.2]{IF} for hypersurfaces in weighted projective spaces,
using the Jacobian ideal $\J$ of the defining polynomial $f_8$ in the
ring of polynomials $\CX[z_0,\ldots,z_5]$ graded according to the weights
of $z_i$ in $\CPW$. It suffices to compute $h^{3,1}(V)=35$.

We obtain from Propositions \ref{betti} and~\ref{sig} that the
asymptotically cylindrical 8-manifold $M_2$ with holonomy $\Spin(7)$ in
this example has $b^4(M_2)=455$ and $b^4_0(M_2)=304$ and $b^4_-(M_2)=72$.
The moduli space for asymptotically cylindrical $\Spin(7)$-structures
on $M_2$ has dimension~224. The manifold $M_2$ is topologically distinct
from the $\Spin(7)$-manifold $M_1$ in~\S\ref{wps}.

The $\Spin(7)$-manifold $M_2$ can also be constructed from
$\CP^4_{1,1,1,1,4}$ with $D$ as in \S\ref{wps} if we choose a second
anticanonical divisor to be $D''=\{z_4^2=0\}\subset\CP^4_{1,1,1,1,4}$ a
linear cone of multiplicity~2. Then the self-intersection of $D$ is
represented by $D\cdot D''=2\Sigma_8$, the Fermat surface with
multiplicity~2. This admissible configuration is
$(\CP^4_{1,1,1,1,4},D,2\Sigma_8,\rho_1)$ a modification of the example
in \S\ref{wps}. In the present case, in order to obtain a 4-orbifold $\tV$
with an anticanonical divisor $\tD$ having holomorphically trivial normal
bundle, {\em two} successive blow-ups of $\Sigma_8$ in $C$ are required as
in Remark~\ref{blow.div}. (Respectively, the Betti numbers of~$M_2$ are
recovered by applying Proposition~\ref{betti} twice.) The first blow-up
produces a 4-orbifold isomorphic to the hypersurface $V$
in~\eqref{hsurf}.


\begin{thebibliography}{CHNP}

\bibitem[APS]{aps}
M.F. Atiyah, V.K. Patodi, and I.M. Singer.
Spectral asymmetry and Riemannian geometry,~I.
{\it Math. Proc. Camb. Phil. Soc.} {\bf 77} (1975), 43--69.

\bibitem[Ba]{baily} W.L. Baily.
The decomposition theorem for $V$-manifolds.
{\it Am. J. Math.} {\bf 78}, (1956) 862--888. 

\bibitem[BM]{biquard-minerbe}
O.~Biquard and V.~Minerbe.
A Kummer construction for gravitational instantons.
{\it Comm. Math. Phys.} {\bf 308} (2011), 773--794. 

\bibitem[Bo]{bonan}
E. Bonan. Sur des vari\'e t\'e s riemanniennes \`a groupe d'holonomie $G_2$
ou $\Spin(7)$. {\it C. R. Acad. Sci. Paris Sér. A} {\bf 262} (1966), 127--129.

\bibitem[Br]{bryant}
R.L.~Bryant.
Metrics with exceptional holonomy.
{\it Ann. of Math.} {\bf 126} (1987), 525--576.

\bibitem[BS]{bryant-salamon}
R.L. Bryant and S.M. Salamon.
On the construction of some complete metrics with exceptional
holonomy. {\it Duke Math. J.} {\bf 58} (1989), 829--850.

\bibitem[CLS]{candelas}
P.~Candelas, M.~Lynker and R.~Schimmrigk.
Calabi--Yau manifolds in weighted $\mathbb{P}_4$.
{\it Nuclear Phys. B} {\bf 341} (1990), 383--402. 

\bibitem[CG]{split.thm}
J.~Cheeger and D.~Gromoll.
The splitting theorem for manifolds of non-negative {R}icci
 curvature. {\it J. Diff. Geom.} {\bf 6} (1971), 119--128.

\bibitem[CHNP]{CHNP}
A. Corti, M. Haskins, J. Nordstr\"om and T. Pacini.
Asymptotically cylindrical Calabi--Yau 3-folds from weak Fano 3-folds.
{\it Geom. Topol.}, to appear.

\bibitem[dR]{deRham}
G. de Rham.
Sur la reductibilité d'un espace de Riemann.
{\it Comment. Math. Helv.} {\bf 26} (1952), 328--344.

\bibitem[Di]{dimca-book}
A.~Dimca. Singularities and topology of hypersurfaces. Springer, 1992.

\bibitem[DD]{dimca85}
A.~Dimca and S.~Dimiev. On analytic coverings of weighted projective
spaces, {\it Bull. London Math. Soc.} {\bf 17} (1985), 234--238.

\bibitem[Do]{idolga}
I.~Dolgachev. Weighted projective varieties.
Group actions and vector fields (Vancouver, B.C., 1981), 34--71.
Lecture Notes in Math. 956. Springer, 1982.

\bibitem[Fe]{fernandez}
M.~Fern{\'a}ndez.
A classification of Riemannian manifolds with structure group $\Spin(7)$.
{\it Ann. Mat. Pura Appl. (4)} {\bf 143} (1986), 101--122.

\bibitem[GPP]{GPP}
G.W.~Gibbons, D.N.~Page and C.N.~Pope.
Einstein metrics on $S^3,\;{\bf {R}}^3$ and ${\bf {R}}^4$ bundles,
{\it Comm. Math. Phys.} {\bf 127} (1990), 529--553.

\bibitem[GH]{GH}
P. Griffiths and J. Harris.
Principles of algebraic geometry.
John Wiley \&\ Sons, 1978.

\bibitem[HL]{HL}
R.~Harvey and  H.B.~Lawson.
Calibrated geometries. {\it Acta Math.} {\bf 148} (1982), 47--157.

\bibitem[HHN]{HHN}
M.~Haskins, H.-J.~Hein, J.~Nordstr\"om.
Asymptotically cylindrical Calabi--Yau manifolds.
arXiv math/1212.6929.

\bibitem[He]{hein}
H.-J.~Hein.
Gravitational instantons from rational elliptic surfaces.
{\it J. Amer. Math. Soc.} {\bf 25} (2012), 355--393.

\bibitem[Ia]{IF}
A.R.~Iano--Fletcher. Working with weighted complete intersections.
Explicit birational geometry of 3-folds, 101--173.
Cambridge Univ. Press, 2000.

\bibitem[Jo1]{joyce-inv}
D.D.~Joyce. Compact 8-manifolds with holonomy Spin(7). {\it Inventiones
Mathematicae} {\bf 123} (1996), 507--552. 

\bibitem[Jo2]{joyce-jdg}
D.D. Joyce. A new construction of compact 8-manifolds with holonomy Spin(7).
{\it J.\ Diff.\ Geom.} {\bf 53} (1999), 89--130.

\bibitem[Jo3]{joyce}
D.D.~Joyce, Compact manifolds with special holonomy. OUP Mathematical
Monographs series, Oxford, 2000.

\bibitem[Ko1]{g2paper}
A.~Kovalev. Twisted connected sums and special Riemannian holonomy.
{\it J. Reine Angew. Math.} {\bf 565} (2003), 125--160.

\bibitem[Ko2]{gokova}
A.~Kovalev. Ricci-flat deformations of asymptotically cylindrical Calabi--Yau
manifolds. Proceedings of G\"okova Geometry-Topology
Conference 2005, 140--156. G\"okova, 2006.

\bibitem[Ko3]{II}
A. Kovalev. Asymptotically cylindrical manifolds with
holonomy $\Spin(7)$. II. In \mbox{preparation}.

\bibitem[KL]{KL}
A.~Kovalev and N.-H.~Lee. $K3$ surfaces with non-symplectic involution and
compact irreducible $G_2$-manifolds.
{\it Math. Proc. Camb. Phil. Soc.}

\bibitem[KN]{KN}
A.~Kovalev and J.~Nordstr\"om. Asymptotically cylindrical 7-manifolds of
holonomy $G_2$ with applications to compact irreducible $G_2$-manifolds.
{\it Ann.\ Glob.\ Anal.\ Geom.} {\bf 38} (2010), 221--257.

\bibitem[Lo]{lockhart}
R.B. Lockhart.
Fredholm, Hodge and Liouville theorems on noncompact manifolds.
{\it Trans. Amer. Math. Soc.}, {\bf 301} (1987), 1--35.

\bibitem[LM]{LM}
R.B. Lockhart and R.C. McOwen.
Elliptic differential operators on noncompact manifolds.
{\it Ann. Scuola Norm. Sup. Pisa Cl. Sci.}, {\bf 12} (1985), 409--447.

\bibitem[MP]{MP}
V.G. Maz'ja and B.A. Plamenevski.
Estimates on $L^p$ and H\"older classes on the
Miranda--Agmon maximum principle for solutions of elliptic boundary value
problems with singular points on the boundary.
{\it Math. Nachr.} {\bf 81} (1978), 25--82. (in Russian)
English translation: Amer. Math. Soc. Transl. Ser. 2 {\bf 123} (1984), 1--56.

\bibitem[Me]{melrose}
R.B. Melrose.
The Atiyah--Patodi--Singer index theorem.
A~K~Peters Ltd., Wellesley, MA, 1993.

\bibitem[MMR]{MMR}
J. Morgan, T. Mrowka and D. Ruberman.
The $L^2$-moduli space and a vanishing theorem for Donaldson polynomial
invariants. International Press, 1994.

\bibitem[Mo]{morrey}
C.B. Morrey.
Multiple integrals in the calculus of variations.
Grundlehren der matematischen Wissenschaften, volume 130.
Springer-Verlag, Berlin, 1966.

\bibitem[No1]{jn-math.proc}
J. Nordstr\"om. Deformations of asymptotically cylindrical $G_2$-manifolds.
{\it Math. Proc. Camb. Phil. Soc.} {\bf 145} (2008), 311--348.

\bibitem[No2]{jn-thesis}
J. Nordstr\"om. Deformations and gluing of asymptotically cylindrical
manifolds with exceptional holonomy. PhD thesis, Cambridge 2008.

\bibitem[SW]{salamon-walpuski}
D.A.~Salamon and T.~Walpuski. Notes on the octonians. arXiv:1005.2820

\bibitem[S]{salamon}
S.M. Salamon.
Riemannian geometry and holonomy groups.
Pitman Res. Notes in Math. 201. Longman, Harlow, 1989.

\bibitem[Sa]{salur}
S. Salur. Asymptotically cylindrical Ricci-flat manifolds.
{\it Proc. Amer. Math. Soc.} {\bf 134} (2006), 3049--3056.

\bibitem[Sat]{satake}
I.~Satake.
On a generalization of the notion of a manifold.
{\it Proc. Nat. Acad. Sci. U.S.A.} {\bf 42} (1956), 359--363.
\end{thebibliography}
\end{document}